\documentclass{amsart}
\usepackage{amsfonts}
\usepackage{amssymb}
\usepackage{amsmath} 
\usepackage[colorlinks=true, allcolors=blue]{hyperref}
\usepackage{cite}
\usepackage{graphicx}

\usepackage{geometry}
\usepackage{tikz-cd}

\usepackage{chngcntr}

\usepackage{pstricks}
\usepackage{xcolor}
\definecolor{dark-green}{RGB}{14,150,2}
\newcommand{\gpoint}{\color{dark-green}{\circ}}
\newcommand{\rpoint}{\red \bullet}

\makeatletter
\DeclareRobustCommand{\pdot}{\mathbin{\mathpalette\pdot@\relax}}
\newcommand{\pdot@}[2]{%
  \ooalign{%
    $\m@th#1\circ$\cr
    \hidewidth$\m@th#1\cdot$\hidewidth\cr
  }%
}
\makeatother
\newcommand{\sspoint}{\blue \pdot}

\usepackage{stmaryrd} 

\theoremstyle{plain}
\newtheorem{theorem}{Theorem}[section]
\newtheorem{lemma}[theorem]{Lemma}
\newtheorem{proposition}[theorem]{Proposition}
\newtheorem{corollary}[theorem]{Corollary}

\counterwithin{figure}{section}

\theoremstyle{definition}
\newtheorem{example}[theorem]{Example}
\newtheorem{definition}[theorem]{Definition}
\newtheorem{remark}[theorem]{Remark}

\newtheorem{notations}[theorem]{Notations}
\newtheorem{setting}[theorem]{Setting}

\begin{document}

\title[]{Recollements for graded gentle algebras from spherical band objects}
\author{Pierre Bodin}
\date{\today}

\begin{abstract}
In this paper we study the localization of a derived category of a graded gentle algebra by a subcategory generated by a spherical band object. This object corresponds to a simple closed curve under the equivalence between the  perfect derived category of the graded gentle algebra and the partially wrapped Fukaya category of the associated graded marked surface, as established by Haiden, Katzarkov and Kontsevich.

We describe this localization as a recollement that involves the derived category of a new graded algebra given by quiver and relations. This leads us to the introduction of the class of graded pinched gentle algebras, a generalization of graded gentle algebras. We then show that these algebras are in bijection with graded marked surfaces with conical singularities. Moreover, under this correspondence the localization process amounts to the contraction of the closed curve.
\end{abstract}

\maketitle

\setcounter{tocdepth}{1}
\tableofcontents

\section{\textbf{Introduction}} \label{sec1}

\subsection{Context}

In this paper, we describe the partially wrapped Fukaya category of a graded marked surface after taking a quotient by certain spherical objects. Geometrically, such an object corresponds to a closed curve on the surface, and taking the quotient will amount to contracting the curve.

\begin{figure}[h]
    \centering
    \includegraphics[width=0.65\linewidth]{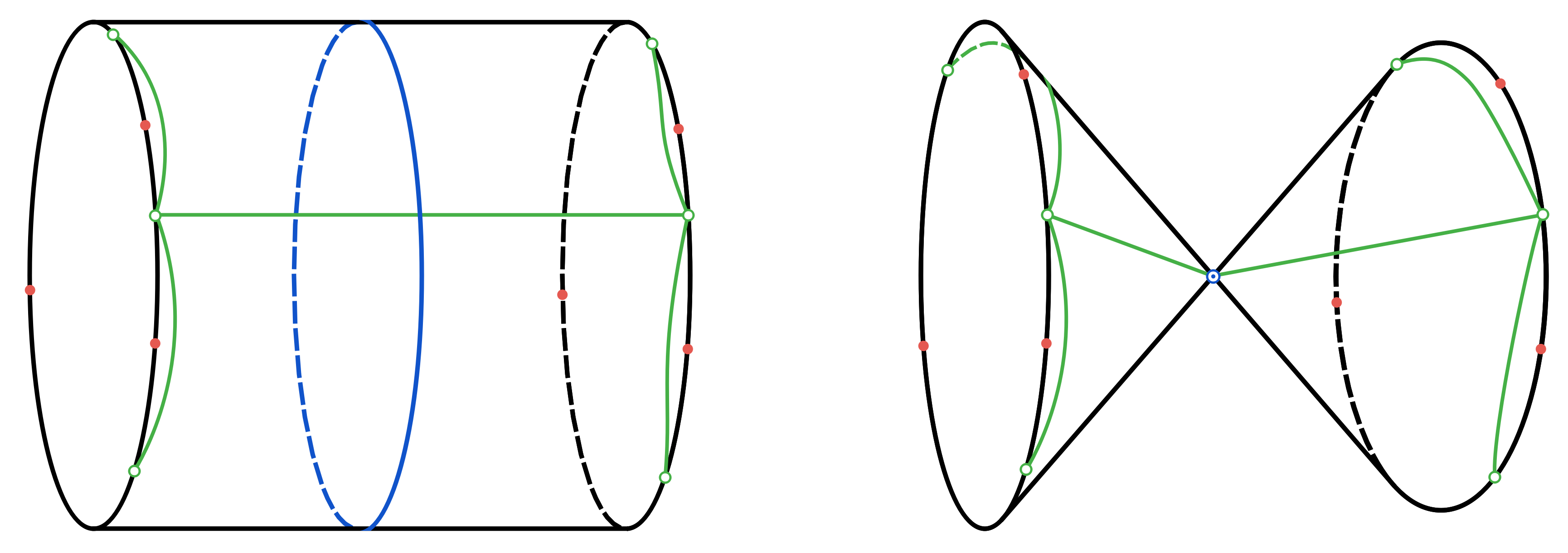}
    \caption{\label{Fig. Ex of surf contraction} A marked surface with a simple closed curve (blue), and its \newline \centering{marked surface with conical singularities obtained by contraction.}}
\end{figure}

In \cite{J22}, the author considers families $f: X \rightarrow \mathbb{C}$ of symplectic manifolds with a singular fiber over~$0$, and defines the wrapped Fukaya category of the singular fiber to be a certain localization of the Fukaya category of a nearby fiber. Our work can be seen as an algebraic analogous construction in which we seek an explicit description of the localization. We will do this through the use of \emph{gentle algebras}.

Since their introduction \cite{AH81,AH82,AS87}, gentle algebras were found to be linked to many other areas of mathematics. Their representation theory has been widely studied, in particular in \cite{WW85,GP68,BR87,CB89,Kra91}. Originating in the theory of cluster algebras from triangulated surfaces \cite{FG09, FST08}, geometric models for gentle algebras have been developed in numerous directions \cite{LF09,ABCJP10,BCS18}. 

Gentle algebras also arose in the context of Homological Miror Symmetry as endomorphism algebras of formal generators in partially wrapped Fukaya categories of graded marked surfaces \cite{HKK17}. Conversely, in works such as \cite{LP20,OPS18,PPP19,BCS18}, any graded gentle algebra is associated to a graded marked surface with graded admissible dissection. In the homologically smooth case, \cite{LP20} showed that the bounded derived category of the algebra in equivalent to the partially wrapped Fukaya category of the associated graded marked surface. Moreover, in the ungraded case, based on the description of the bounded derived category of gentle algebra given in \cite{BM03,BD05,ALP16,CPS19,CPS21}, it is shown in \cite{OPS18} that the marked surface with admissible dissection serves as a model for the indecomposable objects, the irreducible morphisms and their cones, and the Auslander-Reiten triangles.

These geometric models have since been fruitful for the study of the bounded derived categories of gentle algebras. For example, this led to a complete derived invariant for graded gentle algebras \cite{LP20, APS23,Opp19,JSW23}, generalizing a numerical derived invariant of \cite{AAG08}(see also \cite{BH07}).

More generally, many algebraic operations on the derived category find an interpretation in terms of transformations at the level of the surface model. For example, \cite{CS23,CJS23,JSW23} illustrated the usefulness of this correspondence in the study of the silting theory of the derived category. In particular, silting mutation corresponds to the changing of graded arcs and in some cases to the flipping of diagonals in a quadrilateral. Another example is the correspondence established in \cite{Opp19} between spherical twists of the derived category and Dehn twists on the surface.

The study of the localization at a boundary arc was done in \cite{HKK17}, relying on the localization of strictly unital $A_{\infty}$-categories introduced in \cite{LO06}. A more general study of such localization at an admissible collection of arcs was done in \cite{CS23,CJS23}. It is shown that the surface obtained by cutting along the arcs, and what the authors call the surface obtained by supporting along the arcs, give rise to a recollement at the level of the derived categories (see for instance \cite{CJS23}(Definition 1.5) for a definition).

In this paper we study the localization at a simple closed curve. Such a curve corresponds to a spherical object in the Fukaya category. Spherical objects were introduced in \cite{ST01} in order to construct categorical actions of the braid group. Since then, they have been shown to play an important role in the study of triangulated categories, and further similar braid group actions were investigated (see for instance \cite{ST01,KS02,GTW17,AL17,Opp19} for more details).

\vspace{2mm}

We prove that for a gentle algebra $\Lambda$, the localization of the derived category $\mathcal{D}(\Lambda)$ by a spherical band object is equivalent to the derived category of an algebra $\Lambda_{(\alpha, \beta)}$ which we decribe explicitly by quiver and relations (see Definition~\ref{Def. graded pinched gentle}). We call $\Lambda_{(\alpha, \beta)}$ a \emph{pinched gentle algebra}.

It sits in a recollement: (see Theorem~\ref{Th. Main gentle})

\begin{equation*}
    \begin{tikzcd}
    \mathcal{D}(\Lambda_{(\alpha, \beta)}) \arrow[r] & \mathcal{D}(\Lambda) \arrow[l, shift left=0.75ex] \arrow[l, shift right=0.75ex] \arrow[r] & \mathcal{D}(K[x]/(x^2)) \arrow[l, shift left=0.75ex] \arrow[l, shift right=0.75ex]
    \end{tikzcd}.
\end{equation*}

This localization process can be iterated (see Theorems~\ref{Th. Pinched Main recollement} and~\ref{Prop. Pinched Main quotient}). This is reminiscent of the recollements obtained in \cite{CJS23} for localizations with respect to arcs.

On the geometric side, we show in section \ref{Subsec. p graded marked surf} that the isomorphism classes of graded pinched gentle algebras are in bijection with marked surfaces with conical singularities and admissible graded dissections (see Definition~\ref{Def. p graded marked surf} and Proposition~\ref{Prop. 1-1.OPS'Pinched}). This is analogous to  the correspondence established for gentle algebras in \cite{OPS18,BCS18,PPP19}. 

Finally, results on derived equivalences of graded gentle algebras given in \cite{HKK17,LP20} allow us to give a class of derived equivalences between graded pinched gentle algebras (Proposition~\ref{Prop. Easy Pinched Der Eq}). A more in depth analysis of the role of the graded marked surfaces with conical singularities as a model for the bounded derived category of graded pinched algebras will be made in a future work.  

\vspace{2mm}

The paper is structured as follows. In the rest of this section we introduce some definitions and state the main theorems, first for gentle algebras (subsection \ref{Subsec. Main Res Gentle}), then for pinched gentle algebras (subsection \ref{Subsec. Main Res Pinched}). In section \ref{Sec DG} we recall some constructions in DG categories, such as Drinfeld's DG quotient \cite{Dri} which will be used in the proof of the main theorem, and we introduce some notations. Some basic definitions and results on spectral sequences are recalled in section \ref{Sec. Spect Seq}, and applied to describe the first page of a spectral sequence on morphism spaces in DG quotients. In section \ref{Sec. Diss} we first show how one can use the surface model of a gentle algebra in order to choose a set of generators of the derived category which will simplify calculations involving a spherical band object.  Then we introduce the notion of graded marked surfaces with conical singularities and establish the correspondence with graded pinched gentle algebras. The section \ref{Sec. Formality} is devoted to the proof of a technical lemma, namely the formality of the quotient algebra. Finally the proof of the main theorem is made in section \ref{Sec. pf main res}.

\subsection*{Acknowledgements} I am grateful to Giulio Salvatori for sharing his notes on surfacehedra, where the notion of surface with conical singularities is studied. I would also like to thank Raf Bocklandt for pointing me to useful references. This paper is part of my Ph.D. thesis. I would like to thank my advisors Thomas Br\"ustle and Pierre-Guy Plamondon.

\subsection{Definitions and main results}\label{Subsec. Main Res Gentle}

We first introduce some definitions and state our main result in the more restricted context of graded gentle algebras.

\vspace{1mm}

\textbf{Conventions:} A quiver $Q$ is a quadruple $Q = (Q_0, Q_1, \sigma, \tau)$ with finite sets of vertices $Q_0$ and arrows $Q_1$, and with $\sigma, \tau: Q_1 \rightarrow Q_0$ the source and target functions. We consider right modules with the convention of composition (a path $\beta \alpha$ in a path algebra has source $\sigma(\alpha)$ and target $\tau(\beta$)), and we use the cohomological convention for complexes. For each $a \in Q_0$, let $e_a$ be the path of length zero at $a$ and $P_a = e_a \Lambda$ be the associated indecomposable projective $\Lambda$-module. In the rest of the paper $K$ will denote a field.

\begin{definition}\label{Def. Graded Kro and loc at.}
    Let $\Lambda = (KQ / \langle I \rangle, | \ . \ |)$ be a $\mathbb{Z}$-graded algebra.
    \begin{enumerate}
        \item[$\bullet$] $\Lambda$ is said to be \textit{gentle} if:

        \begin{itemize}
            \item[-] Every vertex of $Q$ has at most two incoming and two outgoing arrows,
            \item[-] $I$ is a set of paths of length two satisfying: for all $\alpha \in Q_1$, there is at most one arrow $\beta$ such that $0 \neq \alpha\beta \in I$; at most one arrow $\gamma$ such that $0 \neq \gamma\alpha \in I$; at most one arrow $\beta'$ such that $0 \neq \alpha\beta' \notin I$; at most one arrow $\gamma'$ such that $0 \neq \gamma'\alpha \notin I$.
        \end{itemize}    
    \end{enumerate}

    Suppose now that $\Lambda$ is gentle.
  
    \begin{enumerate}
        \item[$\bullet$] A \textit{graded Kronecker} of $\Lambda$ is a pair of arrows $(\alpha, \beta)$ such that $\alpha$ and $\beta$ have the same source and the same target, the same degree, and such that they are not loops. It is said to be \textit{acyclic} if $\alpha$ and $\beta$ do not belong to an oriented cycle of $(Q, I)$.
    \end{enumerate}

    \begin{enumerate}
        \item[$\bullet$] let $\omega = \alpha + \mu\beta$ for some graded Kronecker $(\alpha, \beta)$ of $\Lambda$ and some $\mu \in K^*$. The \textit{localization of $\Lambda$ at $\omega$} is the graded algebra $\Lambda[\omega^{-1}] = (K\widetilde{Q}/\langle\widetilde{I} \rangle ,\widetilde{| \ . \ |})$ defined by:
            \begin{enumerate}
            \item[-] adding to $Q$ an arrow $\delta$ from 2 to 1 of degree $-|\alpha|$:
            \end{enumerate}
    \end{enumerate}
    
            \begin{center}
            \begin{tikzcd}
            1 \arrow[rr, "\alpha", bend left=25, shift left=2] \arrow[rr, "\beta"', bend right=25, shift right=2] &  & 2 \arrow[ll, "\delta"']
            \end{tikzcd}
            \end{center}

    \begin{enumerate}
        \item[] 
            \begin{enumerate}
            \item[-] adding the relations $\delta \omega - e_1, \omega \delta - e_2$.
            \end{enumerate}
    \end{enumerate}
    
    \begin{enumerate} 
        \item[$\bullet$] Let $(\alpha, \beta)$ be a graded Kronecker of $\Lambda$. Let $1:=\sigma(\alpha)=\sigma(\beta)$, $2:=\tau(\alpha)=\tau(\beta)$, and $\alpha^+, \alpha^-, \beta^+, \beta^-$ be either zero or if they exist, be the (possibly equal) arrows of $Q$ satisfying $\alpha^+ \beta, \beta \alpha^-, \beta^+ \alpha, \alpha \beta^- \in I$ (see Example~\ref{Ex. pinching at a K} left). 
        
        \noindent The \textit{pinching of $\Lambda$ at $(\alpha, \beta)$} is the graded algebra $\Lambda_{(\alpha, \beta)}:= (KQ'/ \langle I' \rangle, | \ . \ |')$ defined by:

        \begin{enumerate}
            \item[-] removing $\alpha$ and $\beta$;
            \item[-] merging 1 and 2 into a new vertex called 1;
            \item[-] adding a loop $\gamma$ of degree 0 at vertex 1;
            \item[-] setting $|\alpha^+|'=|\alpha^+|+|\alpha|$ and $|\beta^+|'=|\beta^+|+|\beta|$, and leaving unchanged the degree of all the other arrows;
            \item[-] keeping all relations in $I$ that don't involve $\alpha$ and $\beta$, and adding the relations  

            \hspace{6,3mm} $\lbrace \alpha^+ \beta^-, \beta^+ \alpha^-, \beta^+(\gamma + e_1), (\gamma + e_1)\beta^-, \alpha^+(\gamma - e_1), (\gamma - e_1)\alpha^- \rbrace$.
        \end{enumerate}
        
        \item[$\bullet$] A \textit{band object supported by $(\alpha, \beta)$} is an element of the perfect derived category $per(\Lambda)$ isomorphic to a shift of a twisted complex (see Notations~\ref{Not. DG recollections}) of the form $(P_2 [|\alpha|] \bigoplus P_1 [1], \partial =
        \begin{pmatrix}
        0 & \alpha + \mu \beta \\
        0 & 0
        \end{pmatrix} )$, for some $\mu \in K^*$. We call $\mu$ the \textit{parameter} of the band object.
    \end{enumerate}
\end{definition}

\begin{example}\label{Ex. pinching at a K} The quiver and relations of a gentle algebra $\Lambda_1$ with a graded Kronecker $(\alpha,\beta)$, and of the associated pinching ${\Lambda_{1}}_{(\alpha,\beta)}$:

\[\begin{tikzcd}
	{\tilde{0}} &&& {\tilde{3}} && {\tilde{3}} && {\tilde{0}} \\
	& 1 & 2 &&&& 1 \\
	0 &&& 3 && 0 && 3 \\
	& {} & {} && {} & {} & {} && {}
	\arrow[""{name=0, anchor=center, inner sep=0}, "{\beta^-}", from=1-1, to=2-2]
	\arrow[""{name=1, anchor=center, inner sep=0}, "{\alpha^-}"', from=3-1, to=2-2]
	\arrow[""{name=2, anchor=center, inner sep=0}, "{\beta^+}", from=2-3, to=1-4]
	\arrow[""{name=3, anchor=center, inner sep=0}, "{\alpha^+}"', from=2-3, to=3-4]
	\arrow[""{name=4, anchor=center, inner sep=0}, "{\alpha^-}", from=3-6, to=2-7]
	\arrow[""{name=5, anchor=center, inner sep=0}, "{\beta^-}", from=1-8, to=2-7]
	\arrow[""{name=6, anchor=center, inner sep=0}, "{\beta^+}", from=2-7, to=1-6]
	\arrow[""{name=7, anchor=center, inner sep=0}, "{\alpha^+}", from=2-7, to=3-8]
	\arrow[""{name=8, anchor=center, inner sep=0}, "\alpha", shift left, from=2-2, to=2-3]
	\arrow[""{name=9, anchor=center, inner sep=0}, "\beta"', shift right, from=2-2, to=2-3]
	\arrow["{\Lambda_1 \ \mathrm{with \ graded \ Kronecker} \ (\alpha,\beta)}"{description}, draw=none, from=4-2, to=4-3]
	\arrow["{{\Lambda_1}_{(\alpha, \beta)}}"{description, pos=1}, draw=none, from=4-5, to=4-6]
	\arrow["{\langle \alpha^+ \beta^-, \ \beta^+ \alpha^-,}"{description, pos=0.4}, shift left=5, draw=none, from=4-6, to=4-9]
	\arrow["{ \ \beta^+(\gamma + e_1), \ (\gamma + e_1)\beta^-,}"{description}, draw=none, from=4-6, to=4-9]
	\arrow["{ \ \alpha^+(\gamma-e_1), \ (\gamma-e_1)\alpha^- \rangle }"{description}, shift right=5, draw=none, from=4-6, to=4-9]
	\arrow[shift right=3, bend left=49, shorten <=5pt, shorten >=5pt, dashed, no head, from=4, to=6]
	\arrow[shift left=3, bend right=49, shorten <=5pt, shorten >=5pt, dashed, no head, from=7, to=5]
	\arrow[shift right=5, bend left=60, shorten <=14pt, shorten >=14pt, dashed, no head, from=0, to=8]
	\arrow[shift left=5, bend right=60, shorten <=14pt, shorten >=14pt, dashed, no head, from=1, to=9]
	\arrow[shift left=5, bend right=60, shorten <=14pt, shorten >=14pt, dashed, no head, from=2, to=8]
	\arrow[shift right=5, bend left=60, shorten <=14pt, shorten >=14pt, dashed, no head, from=3, to=9]
    \arrow["\gamma"', loop, shift right=0.75ex, distance=3em, in=125, out=55,  from=2-7, to=2-7]
\end{tikzcd}\]

\end{example}

Recall that graded gentle algebras are in one-to-one correspondence with graded marked surfaces with admissible dissection, and that under this correspondence, graded curves are associated to objects of the bounded derived category (see subsection \ref{Subsec. Ad adm diss} and \cite{BCS18,OPS18,PPP19} for more details). Our first main result is the following.

\begin{theorem}\label{Th. Main gentle}
    Let $\Lambda'$ be a graded gentle algebra associated to the graded marked surface with admissible dissection $(S,M,\Delta,G)$. Let $\gamma$ be a simple closed curve on $S$ with winding number zero that does not enclose a subsurface containing only punctures, and let $B_{\gamma}$ be an associated band object at the base of a tube.

    \begin{enumerate}
        \item[(1)] There exists a graded gentle algebra $\Lambda$ and an equivalence $\Psi: per(\Lambda') \rightarrow per(\Lambda)$ such that $\Psi(B_{\gamma})$ is a band object supported by an acyclic graded Kronecker $(\alpha, \beta)$ of $\Lambda$, with parameter $\mu$.
        
        \item[(2)] Let $\omega = \alpha + \mu \beta$. There is a recollement:

    \begin{equation*}
    \begin{tikzcd}
    \mathcal{D}(\Lambda[\omega^{-1}]) \arrow[r] & \mathcal{D}(\Lambda) \arrow[l, shift left=0.75ex] \arrow[l, shift right=0.75ex] \arrow[r] & \mathcal{D}(K[x]/(x^2)) \arrow[l, shift left=0.75ex] \arrow[l, shift right=0.75ex]
    \end{tikzcd},
    \end{equation*}

\noindent where $x$ is of degree 1.
        
        \item[(3)] There is an equivalence:

\begin{equation*}
    per(\Lambda')/thick(B_{\gamma}) \simeq per(\Lambda[\omega^{-1}]),
\end{equation*}

\noindent where the left hand side is the Verdier quotient.
    \end{enumerate}

\vspace{1mm}
    
Moreover, if the characteristic of $K$ is different from 2, there is:

    \begin{enumerate}
        \item[(2')] A recollement:

    \begin{equation*}
    \begin{tikzcd}
    \mathcal{D}(\Lambda_{(\alpha, \beta)}) \arrow[r] & \mathcal{D}(\Lambda) \arrow[l, shift left=0.75ex] \arrow[l, shift right=0.75ex] \arrow[r] & \mathcal{D}(K[x]/(x^2)) \arrow[l, shift left=0.75ex] \arrow[l, shift right=0.75ex]
    \end{tikzcd},
    \end{equation*}

\noindent where $x$ is of degree 1.
        
        \item[(3')] An equivalence:

\begin{equation*}
    per(\Lambda')/thick(B_{\gamma}) \simeq per(\Lambda_{(\alpha, \beta)}).
\end{equation*}
        
    \end{enumerate}
\end{theorem}

\begin{remark}
    Theorem~\ref{Th. Main gentle} (1) is a reformulation of Corollary~\ref{Cor. Kro Adapt Gentle Alg}, and (2) and (2') (resp. (3) and (3')) is a direct consequence of Theorem~\ref{Th. Pinched Main recollement} (resp. Theorem~\ref{Prop. Pinched Main quotient}).
\end{remark}

\begin{example}
    The algebras $\Lambda_1$ and ${\Lambda_{1}}_{(\alpha,\beta)}$ of Example~\ref{Ex. pinching at a K} give an example of algebras occurring in the recollement of Theorem~\ref{Th. Main gentle} (2). An illustration of Theorem~\ref{Th. Main gentle} (1) is given by examples \ref{Ex. ad diss} and \ref{Ex. L_0 to L_1} by letting $\Lambda' = \Lambda_0$ and $\Lambda = \Lambda_1$.
\end{example}

\subsection{Graded pinched gentle algebras}\label{Subsec. Main Res Pinched}

The pinching of a graded gentle algebra at a graded Kronecker leads us to the introduction of the class of graded pinched gentle algebras.

\begin{definition}\label{Def. graded pinched gentle}
    A graded $K$-algebra $\Lambda = (KQ/ \langle I\rangle , | \ . \ |)$ is said to be a \textit{graded pinched gentle algebra} if there is a decomposition $Q_1 = Q_1^g \bigsqcup Q_1^p$ and $I = I^g \bigsqcup I^p$ such that 

    \begin{enumerate}
        \item[$\bullet$] $(Q^g:=(Q_0,Q_1^g,\sigma, \tau), \langle I^g \rangle )$ is a gentle bound quiver. The associated gentle algebra is denoted $\Lambda^g:= KQ^g / \langle I^g \rangle$,

        \item[$\bullet$] $Q_1^p$ is a set of degree zero loops supported on different vertices.

        \item[$\bullet$] For $v \in Q_0$, the arrows $\alpha_v^-, \alpha_v^+, \beta_v^-, \beta_v^+ \in Q_1^g$ satisfying $\tau(\alpha_v^-) = \sigma(\alpha_v^+) = v = \sigma(\beta_v^+) = \tau(\beta_v^-)$ and $\beta_v^+ \alpha_v^-, \alpha_v^+ \beta_v^- \in I^g$ (with the possibility of being zero, and with possible compatible identifications between $\lbrace \alpha_v^-, \beta_v^- \rbrace$ and $\lbrace \alpha_v^+, \beta_v^+ \rbrace$) can be chosen such that:
    \end{enumerate}
    
    \begin{center}
        $I^p = \lbrace \beta_v^+(\gamma_v + e_v), (\gamma_v + e_v)\beta_v^-, \alpha_v^+(\gamma_v - e_v), (\gamma_v - e_v)\alpha_v^- \ | \ v \in Q_0, \ \gamma_v \in Q_1^p \ \mathrm{such \ that} \ \sigma(\gamma_v) = v \rbrace$.
    \end{center}
    
\end{definition}

\begin{remark}\label{Rem. Pinched to Gentle}
    The ``pinched" or ``vanishing" relations of $I^p$ are such that  $\forall i \neq j \in Q_0$, $e_j \Lambda e_i \simeq e_j \Lambda^g e_i$.

    In the special case where $\alpha_v^-$ and $\alpha_v^+$ (or equivalently $\beta_v^-$ and $\beta_v^+$) are both zero, the pinched relations at $v$ can be seen as gentle by letting $\gamma_v' = \gamma_v + e_v$.
\end{remark}

\begin{example} The algebra ${\Lambda_1}_{(\alpha,\beta)}$ of Example~\ref{Ex. pinching at a K} is a pinched gentle algebra whose associated gentle algebra $({\Lambda_1}_{(\alpha,\beta)})^g$ is:

    \vspace{-5mm}

\[\begin{tikzcd}
	{\tilde{3}} && {\tilde{0}} \\
	& 1 \\
	0 && 3 \\
	{} & {} & {}
	\arrow[""{name=0, anchor=center, inner sep=0}, "{\alpha^-}", from=3-1, to=2-2]
	\arrow[""{name=1, anchor=center, inner sep=0}, "{\beta^-}", from=1-3, to=2-2]
	\arrow[""{name=2, anchor=center, inner sep=0}, "{\beta^+}", from=2-2, to=1-1]
	\arrow[""{name=3, anchor=center, inner sep=0}, "{\alpha^+}", from=2-2, to=3-3]
	\arrow["{({\Lambda_1}_{(\alpha,\beta)})^g}", draw=none, from=4-1, to=4-3]
	\arrow[shift right=3, bend left=60, shorten <=5pt, shorten >=5pt, dashed, no head, from=0, to=2]
	\arrow[shift left=3, bend right=60, shorten <=5pt, shorten >=5pt, dashed, no head, from=3, to=1]
\end{tikzcd}\]
    
\end{example}

The notion of graded Kronecker and pinching extend naturally to the context of graded pinched gentle algebras.

\begin{definition}\label{Def. Pinched Graded Kro and loc at.}
    Let $\Lambda = (KQ / \langle I \rangle , | \ . \ |)$ be a graded pinched gentle algebra.
    \begin{enumerate}

        \item[$\bullet$] A \textit{graded Kronecker} $(\alpha, \beta)$ of $\Lambda$ is a graded Kronecker of the associated gentle algebra $\Lambda^g$ such that there is no loop $\gamma \in Q_1^p$ based at the source or target of $\alpha$. It is said to be \textit{acyclic} if $\alpha$ and $\beta$ do not belong to an oriented cycle of $(Q^g, I^g)$.

        \item[$\bullet$] For $(\alpha, \beta)$ a graded Kronecker of $\Lambda$, the \textit{pinching of $\Lambda$ at $(\alpha, \beta)$} is the graded pinched gentle algebra $\Lambda_{(\alpha, \beta)}$ obtained by performing the same local transformation of $(Q,I)$ as in Definition~\ref{Def. Graded Kro and loc at.}.
        
        \item[$\bullet$] The notion of a \textit{band object supported by a graded Kronecker} is as in Definition~\ref{Def. Graded Kro and loc at.}.
        
    \end{enumerate}
\end{definition}

We now state the main results of this paper. Since every graded gentle algebra can be seen as a graded pinched gentle algebra, Theorem~\ref{Th. Main gentle} (2) is a particular instance of Theorem~\ref{Th. Pinched Main recollement}, and Theorem~\ref{Th. Main gentle} (3) a particular instance of Theorem~\ref{Prop. Pinched Main quotient}. Stating the results in this generality also allows us to localize at a collection of disjoint simple closed curves by iterating the process.

\begin{theorem}\label{Prop. Pinched Main quotient} Let $\Lambda$ be a graded pinched gentle algebra and $B$ a band object of $per(\Lambda)$ supported by an acyclic graded Kronecker $(\alpha, \beta)$, with parameter $\mu$. There is an equivalence:

    \begin{equation*}
        per(\Lambda)/thick(B) \simeq per(\Lambda[\omega^{-1}]),
    \end{equation*}
    
    \noindent where the left hand side is the Verdier quotient.
    
    \vspace{1mm}
    
    If the characteristic of $K$ is different from 2, there is an equivalence $per(\Lambda[\omega^{-1}]) \simeq per(\Lambda_{(\alpha, \beta)})$.
    
\end{theorem}

\begin{theorem}\label{Th. Pinched Main recollement} Let $\Lambda$ be a graded pinched gentle algebra with an acyclic graded Kronecker $(\alpha, \beta)$, and let $\omega = \alpha + \mu\beta$ for some $\mu \in K^*$. There exists a recollement:

    \begin{equation*}
    \begin{tikzcd}
    \mathcal{D}(\Lambda[\omega^{-1}]) \arrow[r] & \mathcal{D}(\Lambda) \arrow[l, shift left=0.75ex] \arrow[l, shift right=0.75ex] \arrow[r] & \mathcal{D}(K[x]/(x^2)) \arrow[l, shift left=0.75ex] \arrow[l, shift right=0.75ex]
    \end{tikzcd},
    \end{equation*}

\noindent where $x$ is of degree 1.

\vspace{1mm}

If the characteristic of $K$ is different from 2, there is an equivalence $\mathcal{D}(\Lambda[\omega^{-1}]) \simeq \mathcal{D}(\Lambda_{(\alpha, \beta)})$.

\end{theorem}

The proof of Theorem~\ref{Prop. Pinched Main quotient} and Theorem~\ref{Th. Pinched Main recollement} is made in section \ref{Sec. pf main res}.

\section{\textbf{Recollections on DG categories}}\label{Sec DG}

We introduce notations and recall some properties of DG categories. Our mains references are \cite{CC}, \cite{BK} and \cite{Dri}.

\begin{notations}\label{Not. DG recollections} Let $\mathcal{A}$ be a DG category over a base field $K$.

        \vspace{1mm}

        $\bullet$ We denote by  $H^{*}(\mathcal{A})$ its graded homotopy category and by $H^{0}(\mathcal{A})$ its zeroth  homotopy category.

        \vspace{1mm}

        $\bullet$ We recall that the category $\mathcal{A}^{pre-tr}$ of \textit{one-sided twisted complexes} over $\mathcal{A}$ is defined in \cite{BK} as follow:

        \begin{enumerate}
        \item[-] Objects in $\mathcal{A}^{pre-tr}$ are formal expressions $(\bigoplus\limits_{i=1}^{n} C_i[r_i],\partial)$ where $n \geq 0$, $C_i \in \mathcal{A}$, $r_i \in \mathbb{Z}$, $\partial=(\partial_{ij})$, $\partial_{ij} \in \mathcal{A}(C_j,C_i)[r_i - r_j]$ is homogeneous of degree 1, $\partial_{ij}=0$ for $i \geq j$ and $d_{\textrm{naive}}\partial+\partial^2 = 0$, where $d_{\textrm{naive}}\partial:= (d_{\mathcal{A}} (\partial_{ij}))$.
        
        \vspace{3mm}
        
        \item[-] $f = (f_{ij}) \in \mathcal{A}^{pre-tr}( (\bigoplus\limits_{j=1}^{n} C_j[r_j],\partial) , (\bigoplus\limits_{i=1}^{m} C'_i[r'_i],\partial'))$ verify $f_{ij} \in \mathcal{A}(C_j,C'_i)[r'_i - r_j]$ and the composition is the matrix multiplication.
        
        \vspace{3mm}
        
        \item[-] The differential is defined by $df:= d_{\textrm{naive}}f + \partial'f - (-1)^l f\partial$ if deg$f_{ij} = l$.
        \end{enumerate}

    \vspace{3mm}

    The DG category $\mathcal{A}$ can be seen as a full DG subcategory of $\mathcal{A}^{pre-tr}$. For $f: X \rightarrow Y$ a closed morphism of degree 0 of $\mathcal{A}$, let $Cone(f)$ be the object $(Y \bigoplus X[1], 
    \begin{pmatrix}
    0 & f \\
    0 & 0
    \end{pmatrix} ) \in \mathcal{A}^{pre-tr}$. 

    Let $C=(\bigoplus\limits_{j=1}^{n} C_j[r_j],\partial) , C'=(\bigoplus\limits_{i=1}^{m} C'_i[r'_i],\partial') \in \mathcal{A}^{pre-tr}$. A homogeneous basis $\lbrace b_h^{kl} \rbrace_{1 \leq h \leq n^{kl}}$ of $\mathcal{A}(C_l,C'_k)$ for each $l \in \lbrace 1, ..., n \rbrace$ and $k \in \lbrace 1, ..., m \rbrace$ induces a basis~$\lbrace [b_h^{kl}] \ | \ 1 \leq k \leq m,$ and~$1 \leq l \leq n,$ and $1 \leq h \leq n^{kl} \rbrace$ of $\mathcal{A}^{pre-tr}(C,C')$, where $[b_h^{kl}]_{ij} = b_h^{kl}$ if $i=k$ and $j=l$, and zero otherwise. The degrees in $\mathcal{A}^{pre-tr}$ are then $|[b_h^{kl}]|^{pre-tr} = |b_h^{kl}| + r_l - r'_k$.

    We will extend this bracket notation to linear combinations and sometimes drop them when there is no ambiguity and the context is clear.

    \vspace{1mm}

    $\bullet$ The triangulated category $\mathcal{A}^{tr}$ associated to $\mathcal{A}$ is by definition $H^0 (\mathcal{A}^{pre-tr})$. The distinguished triangles are isomorphic to diagrams $X \rightarrow Y \rightarrow Cone(f) \rightarrow X[1]$ with $f: X \rightarrow Y$ a degree 0 closed morphism of $\mathcal{A}$. Any inclusions $\mathcal{A} \hookrightarrow \mathcal{C} \hookrightarrow \mathcal{A}^{pre-tr}$ induce a triangulated equivalence $\mathcal{A}^{tr} \simeq \mathcal{C}^{tr}$.
\end{notations}

\begin{example}\label{Ex. P Lambda tr simeq per}
    Let $\Lambda:= (KQ / \langle I\rangle , | \ . \ |, d)$ be a DG algebra given by a DG quiver with relations. We view $\Lambda$ as a DG category $\mathcal{P}(\Lambda)$ with objects $\mathcal{P}(\Lambda)_0 = Q_0$ and morphism spaces $\mathcal{P}(\Lambda) (i,j) = e_j \Lambda e_i$, with the induced grading and differential. Using \cite{CC}(Rem 6.2.5 (2)), one can see that the split closure of the triangulated category $\mathcal{P}(\Lambda)^{tr}$ is equivalent to the perfect derived category $per(\Lambda)$.
\end{example}

Drinfeld introduced and studied in \cite{Dri} a notion of DG quotient. 

\begin{definition}\cite{Dri}(§3.1) \label{Def DG Quo} Let $\mathcal{A}$ be a DG category and $\mathcal{B}$ a full subcategory.

The \textit{quotient DG category} $\mathcal{A} / \mathcal{B}$ has the same objects as $\mathcal{A}$. It has a new morphism $\epsilon_{U}: U \rightarrow U$ of degree $-1$ for all object $U$ of $\mathcal{B}$, and all of its compositions with existing morphisms. Formally, one has an isomorphism of vector spaces (but not of complexes),

\[
\mathrm{Hom}_{\mathcal{A} / \mathcal{B}}(X,Y) \leftarrow \bigoplus\limits_{n \in \mathbb{N}} \mathrm{Hom}_{\mathcal{A} / \mathcal{B}}^{\langle n \rangle}(X,Y)
\]

\noindent where $\mathrm{Hom}_{\mathcal{A} / \mathcal{B}}^{\langle n \rangle}(X,Y)$ is the direct sum, over all families $(U_i)_{1 \leq i \leq n}$ of objects in $\mathcal{B}$, of tensor products $\mathrm{Hom}_{\mathcal{A}}(U_n,U_{n+1}) \otimes K[1] \otimes \mathrm{Hom}_{\mathcal{A}}(U_{n-1},U_n) \otimes ... \otimes  K[1] \otimes \mathrm{Hom}_{\mathcal{A}}(U_0,U_1)$ where $U_{n+1} = Y, U_0 = X$ and $K[1]$ is the complex with only $K$ in degree -1.

If $\epsilon$ is the  canonical generator of $K[1]$, the application sends the product $f_n \otimes \epsilon \otimes ... \otimes \epsilon \otimes f_0$ to the composition $f_n \epsilon_{U_n} f_{n-1} ... \epsilon_{U_1} f_0$. The differential is given by the Leibniz rule and $d(\epsilon_U) = id_U$ for all $U$ in $\mathcal{B}$. In particular $\mathrm{Hom}_{\mathcal{A} / \mathcal{B}}^{\langle 0 \rangle}(X,Y)=\mathrm{Hom}_{\mathcal{A}}(X,Y)$.
\end{definition}

This DG quotient is an enhancement of the triangulated quotient:

\begin{theorem} \label{Th Dri Quo} \cite{Dri}(Th 3.4)
There is a triangulated equivalence 

\vspace{1mm}

\begin{center}
$(\mathcal{A}/\mathcal{B})^{tr} \simeq \mathcal{A}^{tr}/\mathcal{B}^{tr}.$
\end{center}

\end{theorem}

\section{\textbf{Spectral sequences and Drinfeld quotients}}\label{Sec. Spect Seq}

In this section, $\mathcal{B}$ will denote a full subcategory of a DG category $\mathcal{A}$. We will present how the $\mathrm{Hom}$ spaces of $\mathcal{A}/\mathcal{B}$ naturally inherit a filtration which gives rise to a spectral sequence that can be used to compute the $\mathrm{Hom}$ spaces of $H^* \mathcal{A}/\mathcal{B}$. Our main reference in this section is McCleary's book \textit{A user's guide to spectral sequences} \cite{McC}.

\vspace{2mm}

Using the notation of Definition~\ref{Def DG Quo}, for all objects $X,Y$ in $\mathcal{A}$, the natural filtration
\begin{equation*}
    0 = F^{-1}C \subset ... \subset F^{p-1}C \subset F^p C \subset ... \subset C:= \mathrm{Hom}_{\mathcal{A} / \mathcal{B}}(X,Y)
\end{equation*}
\noindent given by $F^p C \simeq \bigoplus\limits_{n=0}^{p}\mathrm{Hom}_{\mathcal{A} / \mathcal{B}}^{\langle n \rangle}(X,Y)$ is: 

\begin{enumerate}
    \item[-] \textit{Compatible with the differential}: the differential $d$ of $C$ restricts to a map $d: F^p C \rightarrow F^p C$ for each $p \in \mathbb{Z}$,
    \item[-] \textit{Exhaustive}: $C = \bigcup_s F^s C$,
    \item[-] \textit{Bounded below}: for each $n$, there is a value $s(n)$ with $F^{s(n)} C^n = \lbrace 0 \rbrace$.
\end{enumerate}

These data allow us to construct a spectral sequence:

\begin{theorem} \cite{McC}(Th - 2.6)\label{Th. McC Th.2.6}
    The filtration $F$ of $C$ determines a spectral sequence $\lbrace E_{r}^{*,*}, d_r \rbrace$ starting on page 1 with $d_r$ of bidegree $(-r, r+1)$ and
    
    \vspace{1mm}
    \begin{center}
    $E_{1}^{p,q} \simeq  H^{p+q} (F^{p} C / F^{p-1} C)$.
    \end{center}
    
\end{theorem}
    
This spectral sequence is build recursively starting at $r=0$. Let~$E_0^{p,q} = F^p C^{p+q} / F^{p-1} C^{p+q}$ and $d_{0}^{p,q}: E_0^{p,q} \rightarrow E_0^{p,q+1}$ be the differential induced by the quotient.
    
\vspace{2mm}
    
Suppose that $E_r^{p,q}$ and $d_{r}^{p,q}: E_r^{p,q} \rightarrow E_r^{p-r,q+r+1}$ are given. By definition of a spectral sequence, $E_{r+1}^{p,q} \simeq \mathrm{Ker}(d_{r}^{p,q}) / \mathrm{Im}(d_{r}^{p+r,q-r-1})$. One can then verify that the differential $d$ induces a differential $d_{r+1}^{p,q}: E_{r+1}^{p,q} \rightarrow E_{r+1}^{p-(r+1),q+(r+1)+1}$.
   
\vspace{2mm}

The filtration $F$ also induces a filtration of the homology $H(C, d)$ in the following way:
\begin{center}
$F^p H(C,d) =  \mathrm{Im} \big(H(\iota):  H(F^p C, d) \rightarrow H(C,d)\big)$
\end{center}
\noindent where $\iota$ is the inclusion. Recall the following definition:

\begin{definition} \cite{McC}(Def 2.4)
A spectral sequence $\lbrace E_{r}^{*,*}, d_r \rbrace$ is said to \textit{converge} to a graded $R$-module $H^*$ if there is a filtration $F$ on $H^*$ such that
\begin{equation*}
    E_{\infty}^{p,q} \simeq E^{p,q}(H^*,F),
\end{equation*}
\noindent where $E_{\infty}^{p,q}$ is the limit term of the spectral sequence, and $E^{p,q}(H^*,F)$ the bigraded module given by
\begin{equation*}
    E^{p,q}(H^*,F) = F^p H^{p+q} / F^{p-1} H^{p+q}.
\end{equation*}
\end{definition}

The following theorem applies:

\begin{theorem} \cite{McC}(Th 3.2) The spectral sequence of Theorem~\ref{Th. McC Th.2.6} converges to $H(C,d)$, that is,
    \begin{equation*}
        E_{\infty}^{p,q} \simeq F^p H^{p+q} (C,d) / F^{p-1} H^{p+q} (C,d).
    \end{equation*}    
\end{theorem}

Since we are working over a field, all short exact sequences split and the convergence implies that for all $p$ and $q$, 
\begin{equation}\label{Eq. H* DG quot via Ss}
    H^{p+q} (\mathrm{Hom}_{\mathcal{A} / \mathcal{B}}(X,Y)) \simeq \bigoplus\limits_{k \in \mathbb{Z}} E_{\infty}^{p+k,q-k}.
\end{equation}

Note that $\forall p < 0, \forall q \in \mathbb{Z}$, $E_0^{p,q} \simeq 0$. The following proposition will allow us to compute easily the first terms $E_{1}^{p,q}$.

\begin{proposition} \label{prop E^1}
Let $\mathcal{B}$ be a full subcategory of a DG category $\mathcal{A}$ over a field $K$, and let $X,Y$ be two objects of $\mathcal{A}$. Let $\lbrace E_r^{*,*}, d_r \rbrace$ be the spectral sequence associated to the natural filtration of $C = \mathrm{Hom}_{\mathcal{A} / \mathcal{B}}(X,Y)$. Then

\begin{align*}
E_1^{p,q} \simeq \bigoplus\limits_{(U_i)\in \mathcal{B}} \bigoplus\limits_{\substack{k_1 + ... + k_{p+1} \\ = \ 2p + q}} H^{k_{p+1}}(\mathrm{Hom}_{\mathcal{A}}(U_p,Y)) \otimes K[1] \otimes & H^{k_p} (\mathrm{Hom}_{\mathcal{A}}(U_{p-1},U_p)) \otimes ... \\
 & ... \otimes K[1] \otimes H^{k_1} (\mathrm{Hom}_{\mathcal{A}}(X,U_1)). 
\end{align*}

\end{proposition}

\textit{Proof:}

\vspace{3mm}

By definition, $E_1^{p,q} = H^{p+q} (F^p C / F^{p-1} C)$. Moreover $F^p C / F^{p-1} C$ is isomorphic to 

$\bigoplus\limits_{(U_i)\in \mathcal{B}} \mathrm{Hom}_{\mathcal{A}}(U_p,Y) \otimes K[1] \otimes \mathrm{Hom}_{\mathcal{A}}(U_{p-1},U_p) \otimes ... \otimes K[1] \otimes \mathrm{Hom}_{\mathcal{A}}(X,U_1)$ as a graded vector space. We will show that  $F^p C / F^{p-1} C$ is actually isomorphic to it as a complex, endowed with the differential coming from the tensor product. 

Recall that given complexes $(C^r,d^r),...,(C^1,d^1)$, their tensor product is defined by $(C^r \otimes ... \otimes C^1, d^{\otimes})$ where for an homogeneous element $x_r \otimes ... \otimes x_1$, $d^{\otimes}$ is given by $d^{\otimes} (x_r \otimes ... \otimes x_1) = \sum\limits_{k=1}^{r}(-1)^{\sum\limits_{l=k+1}^{r}deg(x_l)} x_r \otimes ... \otimes d_k (x_k) \otimes ... \otimes x_1$.

On one hand for $f_p \otimes \epsilon_p \otimes f_{p-1} \otimes ... \otimes \epsilon_1 \otimes f_0$ a homogeneous element in $(\bigoplus\limits_{(U_i)\in \mathcal{B}} \mathrm{Hom}_{\mathcal{A}}(U_p,Y) \otimes K[1] \otimes \mathrm{Hom}_{\mathcal{A}}(U_{p-1},U_p) \otimes ... \otimes K[1] \otimes \mathrm{Hom}_{\mathcal{A}}(X,U_1) , d^{\otimes})$,

\begin{align*}
d^{\otimes} (f_p & \otimes \epsilon_p \otimes f_{p-1} \otimes ... \otimes \epsilon_1 \otimes f_0)  \\
& =  \sum\limits_{k=0}^{p}(-1)^{\sum\limits_{l=k+1}^{p}deg(f_l)-(p-k)} f_p \otimes \epsilon_p \otimes ... \otimes d (f_k) \otimes ... \otimes \epsilon_1 \otimes f_0 \\
& +  \sum\limits_{k=1}^{p}(-1)^{\sum\limits_{l=k}^{p}deg(f_l)-(p-k)} f_p \otimes \epsilon_p \otimes ... \otimes d (\epsilon_k) \otimes ... \otimes \epsilon_1 \otimes f_0 \\
& =  \sum\limits_{k=0}^{p}(-1)^{\sum\limits_{l=k+1}^{p}deg(f_l)-(p-k)} f_p \otimes \epsilon_p \otimes ... \otimes d (f_k) \otimes ... \otimes \epsilon_1 \otimes f_0
\end{align*}

since $K[1]$ is endowed with the zero differential.

On the other hand, for a homogeneous element $\overline{f_p \epsilon_p f_{p-1} ... \epsilon_1 f_0}$ of $F^p C / F^{p-1} C$, 

\begin{align*}
d(\overline{f_p \epsilon_p f_{p-1} ... \epsilon_1 f_0}) = & \sum\limits_{k=0}^{p}(-1)^{\sum\limits_{l=k+1}^{p}deg(f_l) - (p-k)} \overline{f_p \epsilon_p ... d(f_k) ... \epsilon_1 f_0} \\
& + \sum\limits_{k=1}^{p}(-1)^{\sum\limits_{l=k}^{p}deg(f_l) - (p-k)} \overline{f_p \epsilon_p ... f_k f_{k-1} ... \epsilon_1 f_0} \\
& = \sum\limits_{k=0}^{p}(-1)^{\sum\limits_{l=k+1}^{p}deg(f_l) - (p-k)} \overline{f_p \epsilon_p ... d(f_k) ... \epsilon_1 f_0}
\end{align*}

since the $f_p \epsilon_p ... f_k f_{k-1} ... \epsilon_1 f_0$ are in $F^{p-1} C$. This shows that the differential induced by the quotient coincides with the one coming from the tensor product, inducing an isomorphism of complex.

\vspace{3mm}

For the moment we have shown an isomorphism of complexes

\begin{align*}
E_1^{p,q} \simeq H^{p+q} ( \bigoplus\limits_{(U_i)\in \mathcal{B}} \mathrm{Hom}_{\mathcal{A}}(U_p,Y) \otimes K[1] \otimes & \mathrm{Hom}_{\mathcal{A}}(U_{p-1},U_p) \otimes ... \\
 & ... \otimes K[1] \otimes \mathrm{Hom}_{\mathcal{A}}(X,U_1), d^{\otimes}). 
\end{align*}

We conclude using the Künneth formula for complexes (\hspace{1sp}\cite{Wei}(Th 3.6.3)). Given two complexes $P$ and $Q$, there exist for all $n$ a short exact sequence

\[
0 \rightarrow \bigoplus\limits_{i+j = n} H^i(P) \otimes H^j (Q) \rightarrow H^n (P \otimes Q) \rightarrow \bigoplus\limits_{i+j = n-1} Tor_1^K (H^i (P), H^j (Q) ) \rightarrow 0.
\]

Since we are working with vector spaces, the torsion groups vanish and there is an isomorphism $\bigoplus\limits_{i+j = n} H^i(P) \otimes H^j (Q) \simeq H^n (P \otimes Q)$. This gives

\begin{align*}
E_1^{p,q} \simeq \bigoplus\limits_{(U_i)\in \mathcal{B}} \bigoplus\limits_{\substack{h_1 + ... + h_{2p+1} \\ = p+q}} H^{h_{2p+1}} (\mathrm{Hom}_{\mathcal{A}}(U_p,Y)) \otimes H^{h_{2p}} (K[1]) \otimes & H^{h_{2p-1}} (\mathrm{Hom}_{\mathcal{A}}(U_{p-1},U_p)) \otimes ... \\
 & ... \otimes H^{h_2}(K[1]) \otimes H^{h_1} (\mathrm{Hom}_{\mathcal{A}}(X,U_1)). 
\end{align*}

Finally since $H^i (K[1])$ is non zero if and only if $i = -1$, this gives the desired formula.

\hfill $\square$

\section{\textbf{Admissible dissections and graded pinched gentle algebras}}\label{Sec. Diss}

\subsection{Adapted admissible dissections}\label{Subsec. Ad adm diss}

In order to simplify computations in $per(\Lambda)$, it will be useful to choose an adapted set of generators. Recall that graded gentle algebras are in one-to-one correspondence with marked surfaces with graded admissible dissection, that is, a quadruple $(S, M, \Delta, G)$ where $S$ is a compact oriented surface with boundary, $M$ a finite set of $\gpoint$ and $\rpoint$ marked points and punctures, $\Delta$ an admissible $\gpoint$-dissection, and $G$ a grading of the minimal oriented intersections. See for example \cite{HKK17,OPS18,APS23,LP20} for the statement of this result, and for more details and examples. Throughout this section, we will use the terminology of \cite{CJS23}(section 1.5). Contrary to the convention adopted in \cite{APS23}, we will draw marked surfaces with admissible dissection in such a way that the orientation of the arrows of the corresponding quiver will be given by rotating clockwise around a $\gpoint$-point.

By a theorem of \cite{LP20} this correspondence gives rise, in the homologically smooth case (that is, in the case without $\rpoint$-punctures), to an equivalence between the derived category of the graded gentle algebra and the partially wrapped Fukaya category $\mathcal{W}(S, M, \eta(\Delta, G))$ of the associated marked surface with graded admissible dissection:

\begin{theorem}\label{Th. LP19 per Fuk eq} \cite{LP20}(Theorem 3.2.2)
    For a homologically smooth graded gentle algebra $\Lambda$ with associated marked surface with graded admissible dissection $(S, M, \Delta,G)$, there is a line field $\eta(\Delta,G)$ on $S$ such that there is an equivalence 
    
    \begin{center}
        $per(\Lambda) \simeq \mathcal{W}(S, M, \eta(\Delta, G))$.
    \end{center}
\end{theorem}

\begin{remark} In case $\Lambda$ is not homologically smooth, following \cite{OZ22}(Definition 3.24, Remark 3.26) we can still embed $per(\Lambda)$ into a partially wrapped Fukaya category $\mathcal{W}(S_s,M_s,\eta(\Delta_s,G_s))$ as follows. If $\Lambda$ is associated to $(S,M,\Delta,G)$, let $S_s$ be the marked surface obtained by replacing each $\rpoint$-puncture by a boundary component containing one $\rpoint$-point and one $\gpoint$-point, and let $M_s$ be the new set of marked points and punctures. The collection $\Delta$ induces an arc system on $S_s$ that can be completed into a full arc system $\Delta_s$ by choosing a new $\gpoint$-arc for each new boundary component. Finally a choice of grading $G_s$ that specializes to $G$ on $\Delta$ induces an equivalence between $per(\Lambda)$ and the full subcategory of $\mathcal{W}(S_s,M_s,\eta(\Delta_s,G_s))$ generated by $\Delta$. Equivalences for $\mathcal{W}(S_s,M_s,\eta(\Delta_s,G_s))$ induce equivalences for $per(\Lambda)$ (see \cite{OZ22}, Lemma 3.28).
\end{remark}

For $\gamma: \mathbb{S}^1 \rightarrow S$ an immersed curve, we denote by $w_{\eta(\Delta, G)}(\gamma)$ its winding number with respect to $\eta(\Delta, G)$ (see \cite{LP20} Definition 1.1.3). 

\vspace{2mm}

Following \cite{OPS18}(Assumption 2.7), we assume that any finite collection of curves is in minimal position, that is, the number of intersections of each pair of (not necessarily distinct) curves in this set is minimal in their respective homotopy class. According to \cite{T08}, it follows from \cite{FHS82} and \cite{N.C01} that, up to homotopy, this assumption is always satisfied. Given two curves $\gamma, \gamma'$, denote by $|\gamma \cap \gamma'|$ this minimal number of intersections.

\vspace{2mm}

For $\gamma: \mathbb{S}^1 \rightarrow S$ a simple closed curve, let $D_{\gamma}: S \rightarrow S$ be the associated Dehn twist along $\gamma$ (see for instance \cite{FM12} for a definition). For $\delta$ a $\gpoint$-arc on $(S,M)$ intersecting a simple closed curve $\gamma$, each (possibly equal) end point of $\delta$ give rise to an oriented intersection between $\delta$ and $D_{\gamma}(\delta)$, denoted $\alpha(\delta)$ and $\beta(\delta)$. Recall that an oriented intersection between to consecutive $\gpoint$-arcs at a $\gpoint$-point is called a \textit{minimal oriented intersection}.

\begin{proposition}\label{Prop. Kronecker Adapted Dissection}
    Let $(S,M,\Delta,G)$ be a graded marked surface and let $\gamma_1,...,\gamma_r$ be a collection of pairwise non-intersecting simple closed curves on $S$ with zero winding number. Suppose moreover that these curves do not enclose a subsurface containing only punctures.

    \vspace{1mm}

    There exists a graded admissible $\gpoint$-dissection $(\Delta',G')$ on $(S,M)$ such that $\eta(\Delta',G') \simeq \eta(\Delta,G)$, and for all $i \in \lbrace 1,...,r \rbrace$ there exists $\delta_i \in \Delta'$ such that:
    \begin{itemize}
        \item The endpoints of $\delta_i$ are not $\gpoint$-punctures,
        \item $|\delta_i \cap \gamma_i| = 1$,
        \item $D_{\gamma_i}(\delta_i) \in \Delta'$,
        \item  For all $\delta \in \Delta' \backslash \lbrace \delta_i, D_{\gamma_i}(\delta_i) \rbrace$, $|\delta \cap \gamma_i | = 0$,
        \item $\alpha(\delta_i)$ and $\beta(\delta_i)$ are minimal and $G'(\alpha(\delta_i)) = G'(\beta(\delta_i))$.
    \end{itemize}    
\end{proposition}

\vspace{0.5mm}

We will call such a dissection $(\Delta',G')$ an \textit{admissible dissection adapted to} the collection $\lbrace \gamma_1, ..., \gamma_r \rbrace$.

\begin{example}\label{Ex. ad diss} Consider the graded gentle algebra $\Lambda_0$ given by:

\vspace{-2mm}

\[\begin{tikzcd}
	& 5 & 4 \\
	0 &&& 3 \\
	& 1 & 2
	\arrow["a"', from=2-1, to=3-2]
	\arrow["b"', from=3-2, to=3-3]
	\arrow[""{name=0, anchor=center, inner sep=0}, "c"', from=3-3, to=2-4]
	\arrow[""{name=1, anchor=center, inner sep=0}, "d"', from=2-4, to=1-3]
	\arrow["f", from=2-1, to=1-2]
	\arrow["e", from=1-2, to=1-3]
	\arrow[shift right=3, bend left=60, shorten <=5pt, shorten >=5pt, dashed, no head, from=0, to=1]
\end{tikzcd}\]

\vspace{2mm}

and $|a| = |d| = |f| = 0$, $|c|=1$ and $|b|=|e| = -1$. The associated graded marked surface $S$ is depicted on the left of the following figure:

\vspace{2mm}

\begin{center}
    \includegraphics[scale=0.21]{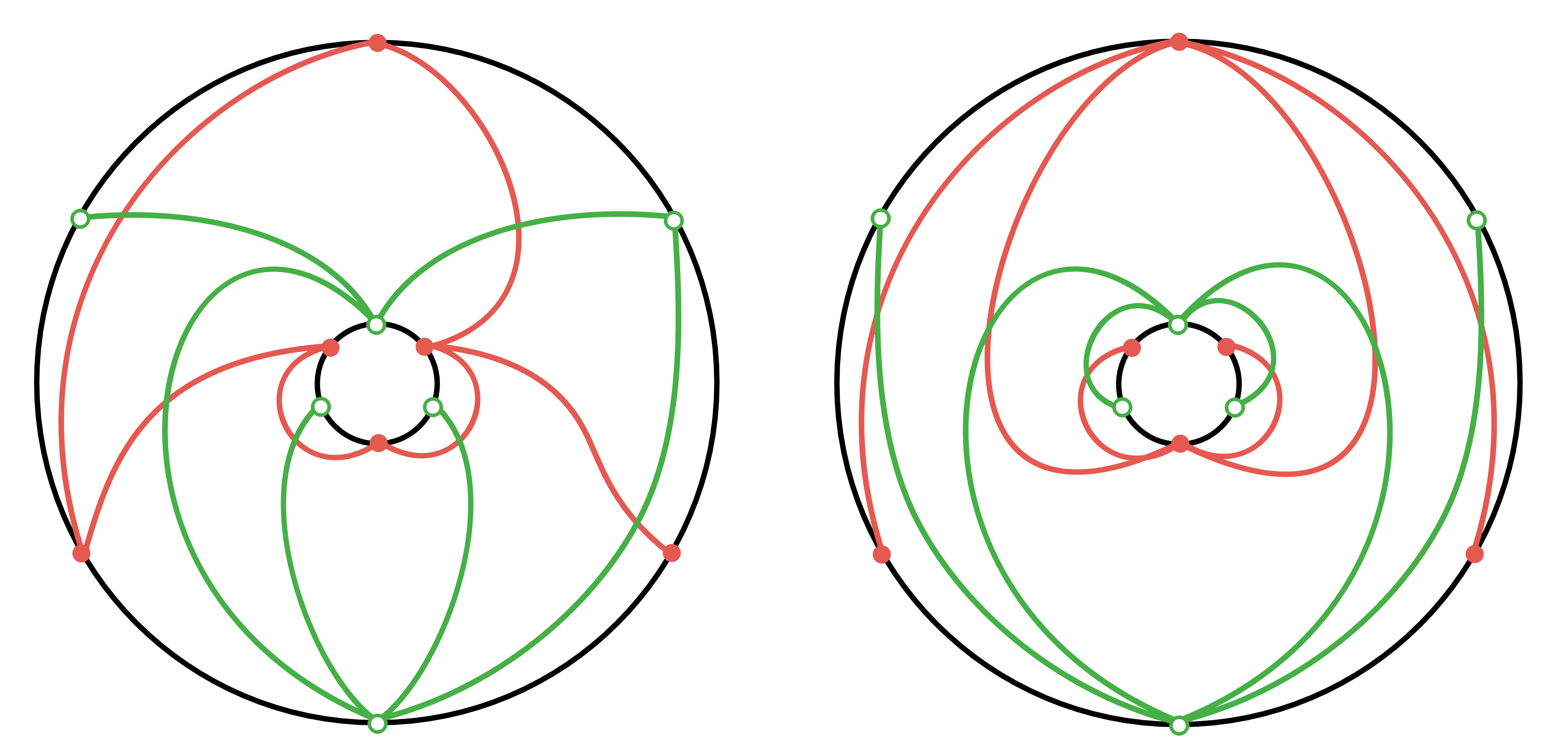}
\end{center}

\vspace{2mm}

The figure on the right is an example of adapted admissible dissection (the grading being zero) to the unique simple closed curve of $S$. Let $S_1$ be this new graded marked surface.

\end{example}

The proof of Proposition~\ref{Prop. Kronecker Adapted Dissection} will rely on the following lemma:

\begin{lemma}\label{Lem. LP.rem1.2.5 winding to genus}\cite{LP20}(Remark 1.2.5) Let $S$ be an oriented surface with non-empty boundary $\partial S$ decomposing as the disjoint union of connected components $\partial S = \bigsqcup\limits_{i=1}^r \partial_i S$. 

A line field $\eta$ on $S$ relates its genus $g$ with the winding number of its boundary components via the formula:

    \begin{center}
        $\sum\limits_{i=1}^r w_{\eta(\Delta,G)} (\partial_i S) = 4 - 2r - 4g$.
    \end{center}
\end{lemma}

\vspace{0.3cm}

\textit{Proof of Proposition~\ref{Prop. Kronecker Adapted Dissection}:} Each simple closed curve $\gamma_i$ has zero winding number, regardless of the choice of orientation.

Let $S'$ be the surface obtained by cutting $S$ along each $\gamma_i$. It is the disjoint union of connected surfaces $S' = \bigsqcup\limits_{j=1}^s S'_j$, and each $\gamma_i$ give rise to two boundary components of $S'$: $\gamma_i^+$ and $\gamma_i^-$. Let us show by contradiction that each $S'_j$ contains a boundary component of the original surface $S$ (and thus contains at least one $\gpoint$ marked point on this boundary).

Since by hypothesis each $S'_j$ cannot contain only punctures, each $S'_j$ containing a puncture must contain a boundary component of $S$. 

Suppose now that $S'_j$ does not contain any puncture and that all its distinct boundary components $\rho_1^j, ..., \rho_{r_j}^j$ are in $\lbrace \gamma_1^+, ..., \gamma_r^+, \gamma_1^-, ..., \gamma_r^- \rbrace$. Let $g_j$ be its genus. By Lemma~\ref{Lem. LP.rem1.2.5 winding to genus}, after choosing the orientation on each $\rho_i^j$ which is compatible with the orientation of $S'_j$, we have 

\begin{align*}
    \sum\limits_{i=1}^{r_j} w_{\eta(\Delta,G)} (\rho_i^j) = 4 - 2r_j - 4g_j & \Leftrightarrow 0 = 4 -2r_j -4g_j \\
    & \Leftrightarrow r_j = 2 \ \text{and} \ g_j = 0,
\end{align*}

but then $\rho_1^j$ and $\rho_2^j$ should be homotopic, a contradiction.

\vspace{0.2cm}

Thus we can choose an arbitrary $\gpoint$ marked point $m_j$ in each $S'_j$.

For all $j \in \lbrace 1, ..., s \rbrace$, and for all $k \in \lbrace 1, ..., r_j \rbrace$, let $\delta_k^j$ be a simple arc from $\delta_k^j (0) = m_j$ to $\delta_k^j (1) = \rho_k^j (1)$ (after the identification with some $\gamma_i^{\pm}: \mathbb{S}^1 \rightarrow S$), and such that $\forall k \neq k' \in \lbrace 1, ..., r_j \rbrace$, $| \delta_k^j \cap \delta_{k'}^j | = 0$.

For each $i \in \lbrace 1, ..., r \rbrace$, this gives rise to a $\gpoint$-arc $\delta_i$ on $S$ intersecting $\gamma_i$ once (and no other $\gamma_{i' \neq i})$ once and transversely, in the following way: if $\gamma_i^+ = \rho_k^j$ and $\gamma_i^- = \rho_{k'}^{j'}$, define $\delta_i$ to be the concatenation $\delta_i = \delta_k^j . (\delta_{k'}^{j'})^{-1}$. Moreover, define $\mu_i^+ = \delta_k^j . \gamma_i . (\delta_k^j)^{-1}$ and $\mu_i^- = \delta_{k'}^{j'} . \gamma_i . (\delta_{k'}^{j'})^{-1}$. Cutting $S$ along the $\mu_i^{\pm}$, we obtain a surface $\Gamma = \Gamma' \bigsqcup ( \bigsqcup\limits_{i = 1}^r \Gamma_i )$, where each $\Gamma_i$ is an annulus containing $\gamma_i$ and whose boundary components $\mu_i^{\pm}$ each have one $\gpoint$ marked point. The following figure depicts the situation:

\vspace{3mm}

\begin{center}
    \includegraphics[scale=0.5]{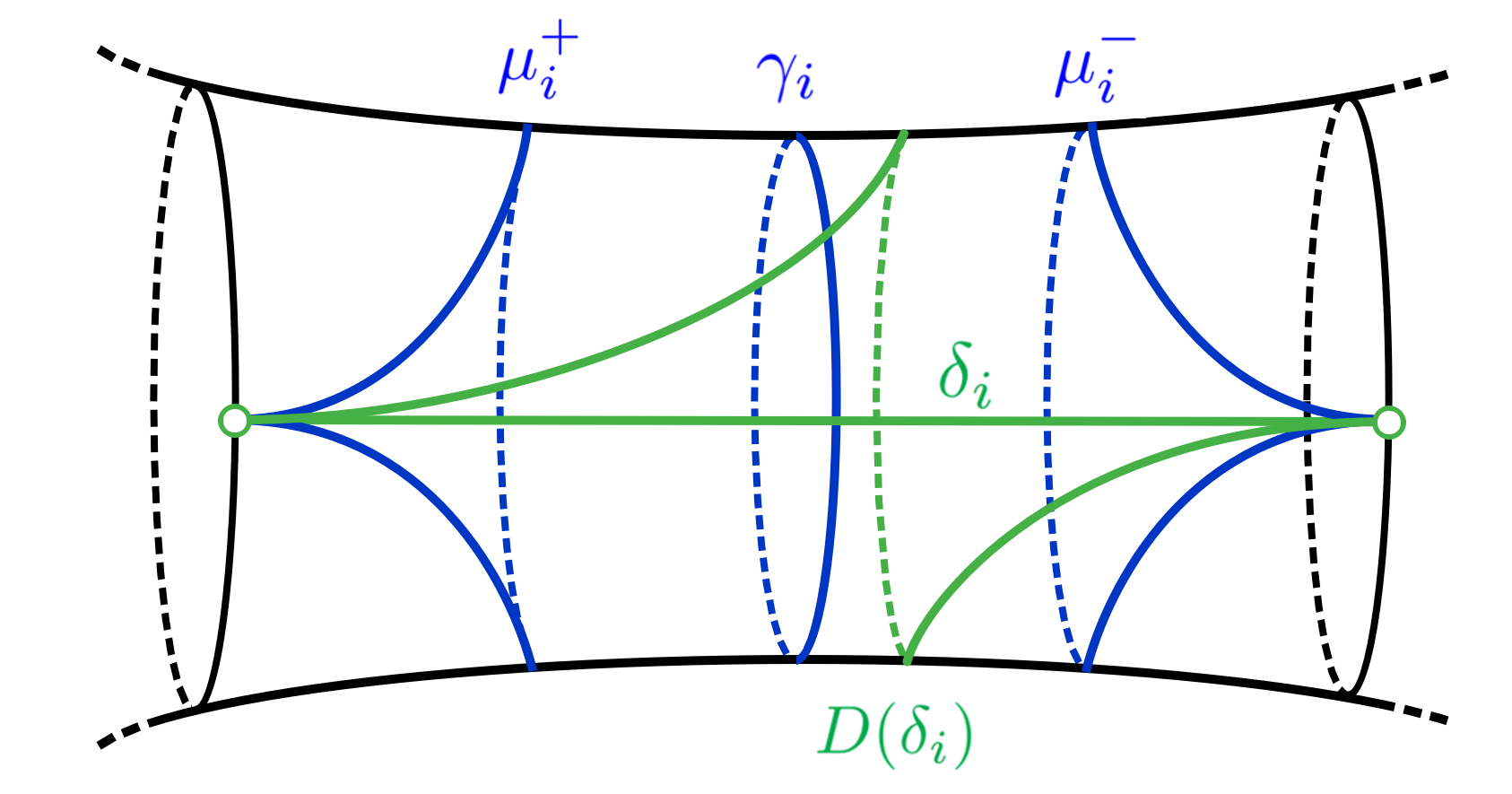}
\end{center}

\vspace{3mm}

Let $\Gamma' = \bigsqcup\limits_{k \in K_0} \Gamma_k$ be the decomposition of $\Gamma'$ into connected components. We have seen above that for each $k$ in $K_0$, $\Gamma_k$ contains at least one $\gpoint$ marked point, thus \cite{HKK17}(Lemma 3.3) ensures that it admits an admissible dissection $A_k$. We set $\Delta':= (\bigsqcup\limits_{k \in K_0} A_k) \bigsqcup \lbrace \delta_1,..., \delta_r, D_{\gamma_1}(\delta_1),..., D_{\gamma_r}(\delta_r) \rbrace$. In order to see that this indeed give an admissible dissection on $(S,M)$, we need to check that each pair of arcs are non-intersecting, which is the case by construction, that they are pairwise distinct and that they cut $S$ into polygons, each of which containing exactly one $\rpoint$-point. We call such a polygon an admissible polygon.

\vspace{3mm}

The curves in $\Delta' \bigsqcup \lbrace \gamma_1,...,\gamma_r \rbrace$ are in minimal position by assumption. For $k \neq k' \in K_0$, any $\delta \in A_k$ and $\delta' \in A_{k'}$ are distinct since they do not share common endpoints. For all $i$, $\delta_i$ is distinct from $D_{\gamma_i}(\delta_i)$ and for $i \neq j \in \lbrace 1,...,r \rbrace$, $\delta_i \not\simeq \delta_j, \delta_i \not\simeq D_{\gamma_j}(\delta_j)$ and $D_{\gamma_i}(\delta_i) \not\simeq D_{\gamma_j}(\delta_j)$ since they do not have the same minimal relative position with respect to the set $\lbrace \gamma_1,..., \gamma_r \rbrace$.

\vspace{3mm}

Let $S_0$ be a subsurface cut out by $\Delta'$. If all boundary arcs of $S_0$ are contained in some $A_k$, then $S_0$ is an admissible polygon since $A_k$ is an admissible collection. Otherwise $S_0$ is obtained by gluing an admissible polygon of some $A_k$ with a collection of triangles of the form $(\mu_i^{\pm}, \delta_i, D_{\gamma_i}(\delta_i) )$. The gluing is made along the $\mu_i^{\mp}$'s, and since elements of $\lbrace \delta_1,..., \delta_r, D_{\gamma_1}(\delta_1),..., D_{\gamma_r}(\delta_r) \rbrace$ are all distinct, $S_0$ is a polygon. Finally since each triangle  $(\mu_i^{\pm}, \delta_i, D_{\gamma_i}(\delta_i) )$ contains no $\rpoint$ points, $S_0$ is admissible.

\vspace{3mm}

Now $\eta(\Delta,G)$ induces a unique grading $G'$ on $\Delta'$ such that $\eta(\Delta',G') \simeq \eta(\Delta,G)$, and by construction $\alpha(\delta_i)$ and $\beta(\delta_i)$ are minimal for all $i$. The fact that $G'(\alpha(\delta_i)) = G'(\beta(\delta_i))$ follows from $w_{\eta(\Delta,G)}(\gamma_i) = 0$.

\hfill $\square$

Proposition~\ref{Prop. Kronecker Adapted Dissection} has the following consequence on the algebraic side:

\begin{corollary}\label{Cor. Kro Adapt Gentle Alg}

    Let $\Lambda = (KQ / I, | \ . \ |)$ be a graded gentle algebra associated to the graded marked surface with admissible dissection $(S,M,\Delta,G)$. 

    Let $B_1,...,B_r$ be a collection of band objects of $per(\Lambda)$ of parameters $(n_i = 1,  \mu_i \in K^*)$ corresponding to a collection $\gamma_1,...,\gamma_r$ of simple closed curves on $S$ under the equivalence $per(\Lambda) \simeq \mathcal{W}(S,M,\eta(\Delta,G))$ or under the inclusion into $\mathcal{W}(S_s,M_s,\eta(\Delta_s,G_s))$ in case $\Lambda$ is not homologically smooth. Suppose that the $\gamma_i$'s are pairwise non-intersecting and do not enclose a subsurface containing only punctures.
    
    \vspace{2mm}
    
    There exists a graded gentle algebra $\Lambda' = (KQ' / I', | \ . \ |')$ and an equivalence $\Psi: per(\Lambda) \rightarrow per(\Lambda')$ such that each $\Psi(B_i)$ is a band object supported by an acyclic graded Kronecker $(\alpha_i,\beta_i)$ of $\Lambda'$, and such that $\forall i \neq j$, $\lbrace \sigma(\alpha_i), \tau(\alpha_i) \rbrace \bigcap \lbrace \sigma(\alpha_j), \tau(\alpha_j) \rbrace = \emptyset$.
\end{corollary}

\begin{example}\label{Ex. L_0 to L_1} Let $\Lambda_0$ be the gentle algebra of Example~\ref{Ex. ad diss}, and $B$ be a band object of parameter $(1,\mu)$ of $per(\Lambda)$ corresponding to the only simple closed curve of $S$. Then Theorem~\ref{Th. LP19 per Fuk eq} shows that $\Lambda_0$ is derived equivalent to the gentle algebra associated the graded marked surface $S_1$. This algebra is the gentle algebra $\Lambda_1$ of Example~\ref{Ex. pinching at a K}.

Under this equivalence, $B$ is isomorphic to (a shift of) $(P_2 \bigoplus P_1 [1], \partial =
        \begin{pmatrix}
        0 & \alpha + \mu \beta \\
        0 & 0
        \end{pmatrix} )$.
    
\end{example}

\textit{Proof of Corollary~\ref{Cor. Kro Adapt Gentle Alg}:} A choice of admissible dissection $(\Delta',G')$ as in \ref{Prop. Kronecker Adapted Dissection} gives the desired algebra $\Lambda'$, which is derived equivalent to $\Lambda$ by \cite{HKK17}(Proposition 3.2). The acyclic graded Kroneckers come from the minimal angles $\lbrace \alpha(\delta_i), \beta(\delta_i) \rbrace$, and the description of the $B_i$'s in $per(\Lambda')$ in term of one sided twisted complexes is given by \cite{HKK17}(Theorem 4.3).

\hfill $\square$

\subsection{Graded marked surfaces with conical singularities}\label{Subsec. p graded marked surf} In this section we establish the correspondence between graded pinched gentle algebras and marked surfaces with conical singularities and simple admissible dissections, a generalisation of graded marked surface to surfaces with conical singularities.

\begin{figure}[h]
\centering
\includegraphics[width=0.91\linewidth]{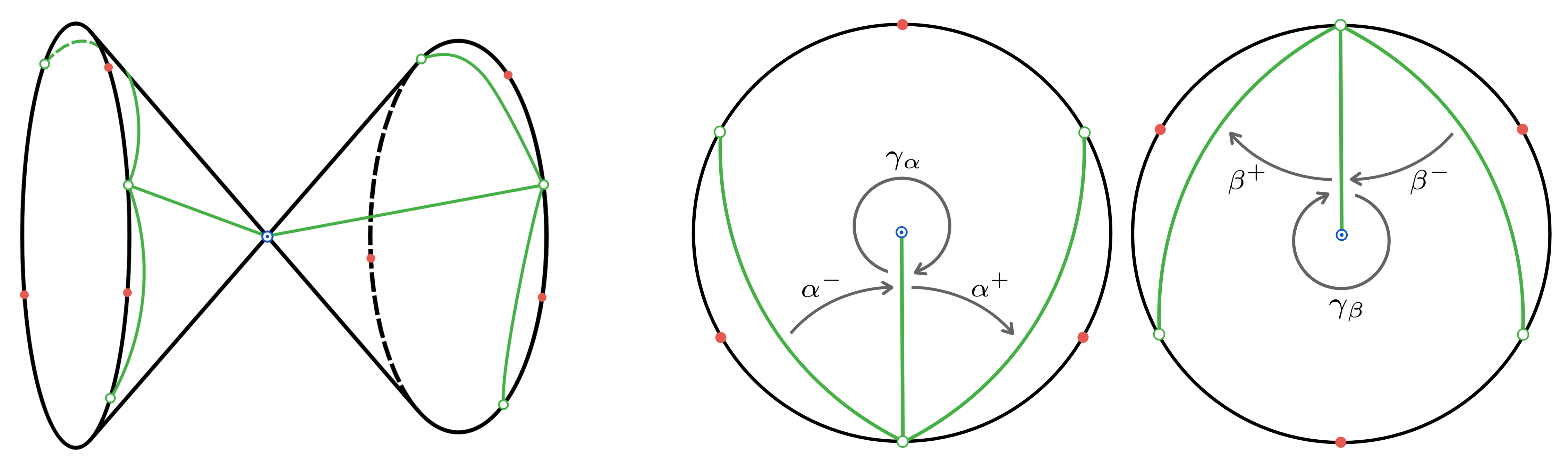}
\caption{\label{Fig. Ex of pgms} A marked surface with conical singularities and graded simple \newline \centering{admissible dissection, and the splitting of its conical singularity.}}
\end{figure}

\begin{definition}\label{Def. p graded marked surf}\
    \begin{enumerate}
        \item[$\bullet$] A \textit{(topological) oriented surface with conical singularities} is a compact topological space $S$ in which every point admits a neighbourhood homeomorphic to one of the following:

        \begin{enumerate}
            \item[-] The open unit disk $D$ of $\mathbb{R}^2$,
        
            \item[-] The half disk $D \cap H$, where $H$ is the closed upper half plane of $\mathbb{R}^2$,

            \item[-] The wedge sum $C_0 = D \sqcup D' / \sim$ of two pointed open disks (see the left-hand side of Figure~\ref{Fig. Ex of pgms}).
        \end{enumerate}

        \noindent Since $S$ is compact, it admits a finite number of boundary components, each homeomorphic to a circle, and a finite set of points $C$, called the set of \textit{conical singularities}, that correspond to the singular point of $C_0$. The surface $S \backslash C$ must be oriented. Elements of $C$ will be represented by symbols $\sspoint$.
        
        \item[$\bullet$] A \textit{marked surface with conical singularities} $\mathcal{S} = (S,C,M \sqcup P)$ is the data of:
        
        \begin{enumerate}
            \item[-] A surface with conical singularities $(S,C)$,
            \item[-] A finite set $P = P_{\gpoint} \sqcup P_{\rpoint}$ of points in the interior of $S \backslash C$, called \textit{punctures} and represented respectively by symbols $\gpoint$ and $\rpoint$,
            \item[-] A finite set $M = M_{\gpoint} \sqcup M_{\rpoint}$ of points in the boundary of $S \backslash C$, called \textit{marked points} and represented respectively by symbols $\gpoint$ and $\rpoint$. They are required to alternate on each boundary component, and each component must contain at least one marked point.
        \end{enumerate}

        \item[$\bullet$] A $\gpoint$-\textit{arc} is a non contractible curve on $S$ with endpoints in $M_{\gpoint} \sqcup P_{\gpoint}$. Every curve will be considered up to homotopy. Note that $\gpoint$-arcs are allowed to go through conical singularities.
        
        \item[$\bullet$] Let $\mathcal{S}_{\sslash}$ be the marked surface obtained from $\mathcal{S}$ by splitting each conical singularity into two $\gpoint$-punctures (see the right-hand side of Figure~\ref{Fig. Ex of pgms}). A \textit{simple admissible dissection} $\Delta$ on $\mathcal{S}$ is a set of $\gpoint$-arcs that don't intersect (in the interior of $S$) and that satisfy: 

        \begin{enumerate}
            \item[-] For each $c \in C$ there is exactly one arc $\gamma_c$ in $\Delta$ that goes through $c$, it goes through it only once and does not go through other singularities,

            \item[-] The collection $\Delta_{\sslash}$ of arcs on $\mathcal{S}_{\sslash}$ induced from $\Delta$ by splitting each $\gamma_c$ in two, is an admissible dissection.
        \end{enumerate}

        \item[$\bullet$] A \textit{grading} $G$ on a simple admissible dissection $\Delta$ is the data of an integer for each minimal oriented angle of $\Delta$.
    \end{enumerate}
\end{definition}

\begin{remark}\label{Rem. Eq pgSurf def} There are two other equivalent ways of thinking of marked surfaces with conical singularities.

\begin{enumerate}
    \item\label{Rem.item. fx-pt free inv and pgms}  First we can think of a marked surface with conical singularities as a marked surface where pairs of $\gpoint$-punctures have been identified. More precisely, given a marked surface $(S,M \sqcup P)$ and a fixed-point free involution $\iota$ on a subset $P_{\gpoint}^{\iota}$ of $P_{\gpoint}$, the associated marked surface with conical singularities $\mathcal{S}$ is obtained by taking the quotient of $S$ under the identifications $p \sim \iota (p)$ for all $p \in P_{\gpoint}^{\iota}$. Thus the set of conical singularities is the quotient of $P_{\gpoint}^{\iota}$ under the action of $\iota$, and the new set of $\gpoint$-punctures is $P_{\gpoint} \backslash P_{\gpoint}^{\iota}$. This construction is inverse to the splitting of the conical singularities $\mathcal{S}_{\sslash}$, where $\iota$ is given by remembering each pairing of new $\gpoint$-punctures. 
    
    Under this correspondence, simple admissible dissections of the marked surface with conical singularities $\mathcal{S}$ are in bijection with admissible dissections on $(S,M \sqcup P)$ satisfying that each puncture in $P_{\gpoint}^{\iota}$ is of valency one. Under this correspondence, the $\gpoint$-arcs $\gamma_p$ and $\gamma_{\iota(p)}$ admitting respectively $p$ and $\iota(p)$ in $P_{\gpoint}^{\iota}$ as one of their end points, give rise via concatenation to an arc going through the corresponding conical singularity.

    Note that here we must exclude the case where $(S,M \sqcup P)$ is a sphere without boundary, with one $\rpoint$-puncture and two $\gpoint$-punctures $p$ and $p'$, and with $\iota(p)=p'$, since the $\gpoint$-arc in the admissible dissection of $(S,M \sqcup P)$ would descend to a closed curve on $\mathcal{S}$. This is the only case where this can happen.

    Each grading $G$ on $(\mathcal{S},\Delta)$ corresponds to a grading $G_{\sslash}$ on $(\mathcal{S}_{\sslash},\Delta_{\sslash})$ that associate zero to each minimal oriented angle at $\gpoint$-punctures arising from conical singularities.

    \item Similarly each collection of non intersecting simple closed curves on a marked surface $(S,M \sqcup P)$ gives rise to a marked surface with conical singularities by taking the topological quotient that identifies each simple closed curve with a point.

\end{enumerate}
\end{remark}

The next proposition follows easily from the correspondence established in \cite{HKK17,BCS18,PPP19,LP20} for gentle algebras. A quiver with relation $(Q,I)$ as in Definition~\ref{Def. graded pinched gentle} will be called a pinched gentle quiver.

\begin{proposition}\label{Prop. 1-1.OPS'Pinched}
    Graded pinched gentle quivers are in one to one correspondence with marked surfaces with conical singularities endowed with a graded simple admissible dissection.
\end{proposition}

\textit{Proof:} Let $(Q,I,| \ . \ |)$ be a graded pinched gentle quiver with set of vanishing loops $Q_1^p$. Let $\gamma_v$ be a loop in $Q_1^p$ attached at vertex $v$, and let $\alpha^+,\alpha^-,\beta^+,\beta^-$ be the (possibly identified or equal to zero) arrows given in Definition~\ref{Def. graded pinched gentle}. By Remark~\ref{Rem. Pinched to Gentle}, we can assume that at least one of $\alpha^+,\alpha^-$, and at least one of $\beta^+,\beta^-$ is non-zero. Let $(Q',I',| \ . \ |')$ be the graded quiver with relations obtained by resolving each loop in $Q_1^p$ in the following way:

\[\begin{tikzcd}
	{\tilde{3}} && {\tilde{0}} &&& {} && {} & {\tilde{3}} && {\tilde{0}} \\
	& 1 && {} & {} && {1_{\alpha}} &&& {1_{\beta}} \\
	0 && 3 &&& 0 && 3 & {} && {} \\
	{} & {} & {} &&&&& {} & {}
	\arrow[""{name=0, anchor=center, inner sep=0}, "{\alpha^-}", from=3-1, to=2-2]
	\arrow[""{name=1, anchor=center, inner sep=0}, "{\beta^-}", from=1-3, to=2-2]
	\arrow[""{name=2, anchor=center, inner sep=0}, "{\beta^+}", from=2-2, to=1-1]
	\arrow[""{name=3, anchor=center, inner sep=0}, "{\alpha^+}", from=2-2, to=3-3]
	\arrow["{(Q,I,| \ . \ |)}", draw=none, from=4-1, to=4-3]
	\arrow[""{name=4, anchor=center, inner sep=0}, "{\alpha^-}", from=3-6, to=2-7]
    \arrow["", maps to, from=2-4, to=2-5] 
    \arrow[""{name=5, anchor=center, inner sep=0}, "{\alpha^+}", from=2-7, to=3-8]
	\arrow[""{name=6, anchor=center, inner sep=0}, "{\beta^-}", from=1-11, to=2-10]
	\arrow[""{name=7, anchor=center, inner sep=0}, "{\beta^+}", from=2-10, to=1-9]
	\arrow["{(Q',I',| \ . \ |')}", draw=none, from=4-8, to=4-9]
	\arrow[""{name=8, anchor=center, inner sep=0}, draw=none, from=2-7, to=1-6]
	\arrow[""{name=9, anchor=center, inner sep=0}, draw=none, from=1-8, to=2-7]
	\arrow[""{name=10, anchor=center, inner sep=0}, draw=none, from=3-9, to=2-10]
	\arrow[""{name=11, anchor=center, inner sep=0}, draw=none, from=2-10, to=3-11]
	\arrow[shift right=3, bend left=49, shorten <=5pt, shorten >=5pt, dashed, no head, from=0, to=2]
	\arrow[shift left=3, bend right=49, shorten <=5pt, shorten >=5pt, dashed, no head, from=3, to=1]
    \arrow[shift left=5, bend right=49, shorten <=8.5pt, shorten >=2pt, dashed, no head, from=5, to=9]
	\arrow[shift right=5, bend left=49, shorten <=8.5pt, shorten >=2pt, dashed, no head, from=4, to=8]
    \arrow[shift right=5, bend left=49, shorten <=2pt, shorten >=8.5pt, dashed, no head, from=10, to=7]
	\arrow[shift left=5, bend right=49, shorten <=2pt, shorten >=8.5pt, dashed, no head, from=11, to=6]
    \arrow["\gamma"', loop, shift right=0.75ex, distance=3em, in=125, out=55,  from=2-2, to=2-2]
    \arrow["\gamma_{\alpha}", loop, shift right=0.50ex, distance=4em, in=55, out=125,  from=2-7, to=2-7]
    \arrow["\gamma_{\beta}", loop, shift right=0.75ex, distance=4em, in=235, out=305,  from=2-10, to=2-10]
\end{tikzcd}\]

Each loop $\gamma$ in $Q_1^p$ gives rise to two degree zero loops $\gamma_{\alpha}$ and $\gamma_{\beta}$ attached on different vertices, and satisfying the relations $\alpha^+ \gamma_{\alpha} = \gamma_{\alpha} \alpha^- = 0 = \beta^+ \gamma_{\beta} = \gamma_{\beta} \beta^-$. Thus $(Q',I',| \ . \ |')$ is a graded gentle quiver and thus corresponds to a marked surface with graded admissible dissection $(S,M \sqcup P,\Delta,G)$. Moreover for each $\gamma$ in $Q_1^p$, the set of couples $(\gamma_{\alpha},\gamma_{\beta})$ give the data of a partial coupling on the $\gpoint$-punctures of $(S,M \sqcup P)$. By Remark~\ref{Rem. Eq pgSurf def} (\ref{Rem.item. fx-pt free inv and pgms}), this gives rise to a unique marked surfaces with conical singularities with graded simple admissible dissection.

Conversely, starting from a marked surfaces $\mathcal{S}$ with conical singularities and graded simple admissible dissection, one can split its conical singularities to obtain $\mathcal{S}_{\sslash}$ as in Remark~\ref{Rem. Eq pgSurf def} (\ref{Rem.item. fx-pt free inv and pgms}). One gets an associated graded gentle quiver and can do the inverse procedure of the last figure to obtain a graded pinched gentle quiver. By definition the two constructions are inverse from each other.

\hfill $\square$

\begin{example}
    Consider the pinched gentle algebra ${\Lambda_{1}}_{(\alpha,\beta)}$ of Example~\ref{Ex. pinching at a K} (it corresponds to the left quiver with relations of the last figure). The corresponding marked surface with conical singularities, depicted on the left part of the following figure, is homeomorphic to the conical singularity $C_0$. The right part illustrates the marked surface with admissible dissection obtained by splitting the conical singularity.  
\end{example}

\begin{remark} Simple admissible dissections are expected to be a special case of a more general notion of admissible dissection on marked surfaces with conical singularities. This notion should induce a bijection with generating sets of the bounded derived category of the corresponding graded pinched gentle algebra. For example, one could consider dissections admitting several arcs going through a conical singularity.
\end{remark}

The localization described in Theorem~\ref{Th. Pinched Main recollement} corresponds at the level of surfaces with conical singularities to the contraction of the simple closed curve. More precisely, we have the following proposition:

\begin{proposition}\label{Prop. pgms comes from curve contra}
    Let $\Lambda$ be a graded pinched gentle algebra with a graded Kronecker $(\alpha,\beta)$, let $\mathcal{S}$ be its associated marked surface with conical singularities and let $\gamma$ be the simple closed curve on $\mathcal{S}$ corresponding to $(\alpha,\beta)$. (Note that by definition of a graded Kronecker, $\gamma$ does not go through a conical singularity).
    
    The marked surface with conical singularities $\mathcal{S}_{(\alpha,\beta)}$, associated to the pinching $\Lambda_{(\alpha,\beta)}$ of $\Lambda$ at $(\alpha,\beta)$, is homeomorphic to the quotient of $\mathcal{S}$ that identifies all points in $\gamma$.    
\end{proposition}

\textit{Proof:} The marked surface with conical singularities obtained by identifying all points in $\gamma$ is the same as the one obtained by:

\vspace{1mm}

(i) Splitting $\mathcal{S}$ at $\gamma$ to obtain a marked surface with conical singularities with new boundary components $\gamma_{\alpha}$ and $\gamma_{\beta}$, 

(ii) contracting $\gamma_{\alpha}$ and $\gamma_{\beta}$ to obtain two $\gpoint$-punctures $p_{\alpha}$ and $p_{\beta}$,

(iii) gluing $p_{\alpha}$ and $p_{\beta}$ to create a conical singularity.

\vspace{1mm}

Splitting the other conical singularities of the marked surface with conical singularities $\mathcal{S}'$ obtained at step (ii) yields exactly the one used in the proof of Proposition~\ref{Prop. 1-1.OPS'Pinched} to define $\mathcal{S}_{(\alpha,\beta)}$. 

\hfill $\square$

\begin{remark}\label{Rem. pinching dissections}
    Up to derived equivalence we can suppose that $(\alpha,\beta)$ is in degree zero. In this case, the pinching of $\mathcal{S}$ at $\gamma$ sends the graded simple admissible dissection associated to $\Lambda$ to the graded simple admissible dissection associated to $\Lambda_{(\alpha,\beta)}$ (with the identification of the arcs corresponding to $\sigma(\alpha)$ and $\tau(\alpha)$).
\end{remark}

We now state a useful result on derived equivalences for graded pinched gentle algebras. For a closed curve $\gamma$ on a marked surface with conical singularities with simple admissible dissection $(\mathcal{S},\Delta,G)$ that does not go through a singularity, we define its winding number to be $w_{\eta(\Delta_{\sslash},G_{\sslash})}(\gamma)$.

\begin{proposition}\label{Prop. Easy Pinched Der Eq}
    Let $\Lambda$ be a graded pinched gentle algebra corresponding to a marked surface with conical singularities with graded simple admissible dissection $(\mathcal{S},\Delta,G)$. Let $\gamma$ be a simple closed curve on $\mathcal{S}$ with zero winding number that does not go through a singularity nor encloses a subsurface containing only punctures or singularities.
    
    There exists a graded simple admissible dissection $(\Delta',G')$ on $\mathcal{S}$ that is adapted to $\gamma$ (see Proposition~\ref{Prop. Kronecker Adapted Dissection}), and such that $per(\Lambda) \simeq per(\Lambda')$ where $\Lambda'$ is the graded pinched gentle algebra associated to $(\mathcal{S},\Delta',G')$.
\end{proposition}

We can then use this new description of the perfect derived category to apply Theorem~\ref{Prop. Pinched Main quotient} and take the pinching of $\Lambda'$ at the graded Kronecker $(\alpha,\beta)$ corresponding to $\gamma$. By Proposition~\ref{Prop. pgms comes from curve contra}, $\Lambda'_{(\alpha,\beta)}$ gives a graded simple admissible dissection on the marked surface with conical singularities obtained by contracting $\gamma$ in $\mathcal{S}$.

\vspace{1mm}

\textit{Proof:} Let $\Gamma$ be the graded gentle algebra obtained by replacing each vanishing loop $\delta_i$ in $Q_1^p$ by a zero graded Kronecker $(\alpha_i,\beta_i)$ in such a way that $\Gamma$ satisfies $\Gamma_{(\alpha_i,\beta_i)_i}=\Lambda$. By Theorem~\ref{Prop. Pinched Main quotient}, each choice of band objects $B_i$ supported on $(\alpha_i,\beta_i)$ gives an equivalence $per(\Lambda) \simeq per(\Gamma) / thick(B_i)_i$.

Let $(\mathcal{S}(\Gamma), \Delta(\Gamma), G(\Gamma))$ be the marked surface with admissible dissection associated to $\Gamma$, and for all $i$ let $\gamma_i$ be the simple closed curve corresponding to $(\alpha_i,\beta_i)$. The curve $\gamma$ is lifted to a curve on $\mathcal{S}$ with zero winding number, and it does not intersect the $\gamma_i$'s. Moreover this collection of closed curves does not enclose a subsurface containing only punctures, so we can apply Proposition~\ref{Prop. Kronecker Adapted Dissection}. Let $(\Delta(\Gamma)',G(\Gamma)')$ be the resulting graded admissible dissection, let $\Gamma'$ be the associated graded gentle algebra, and let $(\alpha'_i,\beta'_i)_i,(\alpha,\beta)$ be the new graded Kroneckers. 

Setting $\Lambda' = \Gamma'_{(\alpha'_i,\beta'_i)_i}$ we get $per(\Lambda') \simeq per(\Gamma')/ thick(B_i)_i \simeq per(\Gamma)/ thick(B_i)_i \simeq per(\Lambda)$. By Remark~\ref{Rem. pinching dissections}, the graded simple admissible dissection associated to $\Lambda'$ corresponds to the image of $(\Delta(\Gamma)',G(\Gamma)')$ under the pinching at each $\gamma_i$, and thus is adapted to $\gamma$. 

\hfill $\square$

\section{\textbf{Formality of the quotient}}\label{Sec. Formality}

In this section we place ourselves in the following setting:
\begin{setting}\label{Set. mc(A)} 
    Let $\Lambda = (KQ/ \langle I\rangle , | \ . \ |)$ be a graded pinched gentle algebra (see Definition~\ref{Def. graded pinched gentle}). Let $B$ be a band object of $\mathcal{P}(\Lambda)^{pre-tr}$ (see Example~\ref{Ex. P Lambda tr simeq per}), of parameter $\mu \in K^*$ and supported by an acyclic graded Kronecker $(\alpha, \beta)$ of $\Lambda$ (see Definition~\ref{Def. Pinched Graded Kro and loc at.}). We set $\omega:= \alpha + \mu \beta$ and $1:=\sigma(\alpha), 2:=\tau(\alpha)$.
    
    Let $\mathcal{A}$ be the full subcategory of $\mathcal{P}(\Lambda)^{pre-tr}$ supported on $\mathcal{P}(\Lambda)_0$ and $B$. The full subcategory $\mathcal{B}$ of $\mathcal{A}$ supported on $B$ satisfies the following property:
\end{setting}

\begin{theorem}\cite{BK}(§4 - Th 1)\label{Th. BK thick} Let $\langle B\rangle$ be the smallest strictly full triangulated subcategory of $\mathcal{A}^{tr} \simeq per(\Lambda)$ containing $B$. There is an equivlance:
    \begin{equation*}
        \mathcal{B}^{tr} \simeq \langle B\rangle.
    \end{equation*}
\end{theorem}

We will show that $(\mathcal{A} / \mathcal{B} )_{\circ}$, the full subcategory of $\mathcal{A} / \mathcal{B}$ supported on $\mathcal{P}(\Lambda)_0$, is formal.

\subsection{A quasi-equivalence}\

\begin{lemma}\label{Lem. A/B circ P(quiver) description}
    There is a DG equivalence

    \begin{center}
        $(\mathcal{A} / \mathcal{B} )_{\circ} \simeq \mathcal{P}(K\widehat{Q}/\langle I\rangle , \widehat{| \ . \ |}, \widehat{d})$,
    \end{center}

    where:
    \begin{enumerate}

        \item[-] $(K\widehat{Q}/\langle I\rangle , \widehat{| \ . \ |})$ is obtained by adding to $Q$ an arrow $\gamma$ from 1 to 2 of degree $|\alpha|-2$, loops $\epsilon_1$ and $\epsilon_2$ at 1 and 2 of degree $-1$, and an arrow $\delta$ from 2 to 1 of degree $-|\alpha|$:
    \end{enumerate}
    
        \begin{center}
        \begin{tikzcd}
        1 \arrow["\epsilon_1", loop, distance=2em, in=145, out=215] \arrow[rr, "\alpha", shift left] \arrow[rr, "\beta"', shift right] \arrow[rr, "\gamma", bend left=35, shift left=2] &  & 2 \arrow["\epsilon_2", loop, distance=2em, in=325, out=35] \arrow[ll, "\delta", bend left=35, shift left=2]
        \end{tikzcd}
        \end{center}

    \begin{enumerate}
        \item[-] $\widehat{d}$ is defined by setting $\widehat{d}(\epsilon_1)=\delta \omega - e_1$, $\widehat{d}(\epsilon_2)= (-1)^q(e_2 - \omega \delta)$, $\widehat{d}(\gamma)=(-1)^{q+1}(\omega\epsilon_1+\epsilon_2\omega)$, and by sending every other arrows to zero.
    \end{enumerate}
\end{lemma}

\textit{Proof:} Recall that $\epsilon$ is the new generator added to $\mathrm{End}_{\mathcal{A}}(B)$ when taking the quotient.

Let $i,j \in Q_0$. Using the bracket of Notation~\ref{Not. DG recollections}, let $[e_1):=[e_1] \in \mathrm{Hom}_{\mathcal{A}}(P_1,B)$, $[e_2):=[e_2] \in \mathrm{Hom}_{\mathcal{A}}(P_2,B)$, $(e_1]:=[e_1] \in \mathrm{Hom}_{\mathcal{A}}(B,P_1)$ and $(e_2]:=[e_2]\in\mathrm{Hom}_{\mathcal{A}}(B,P_2)$. Every morphism $\rho = 
\begin{pmatrix}
    \rho_2 \\
    \rho_1
\end{pmatrix} \in \mathrm{Hom}_{\mathcal{A}}(P_i,B)$ factors as $\rho = (e_2]\rho_2 + (e_1]\rho_1$ and dually for $\rho' \in \mathrm{Hom}_{\mathcal{A}}(B,P_j)$.

The equivalence is then given by sending $\mathrm{Hom}_{\Lambda}(P_i,P_j)$ to itself via $\Lambda \hookrightarrow K\widehat{Q}/\langle I \rangle $ for all $i,j$, and $(e_1]\epsilon[e_1)$ to $\epsilon_1$, $(e_2]\epsilon[e_2)$ to $\epsilon_2$, $(e_1]\epsilon[e_2)$ to $\delta$ and $(e_2]\epsilon[e_1)$ to $\gamma$.

\hfill $\square$

\begin{definition}\label{Def. loc at gr K and phi to be q.e.} Recall that the localization $\Lambda[\omega^{-1}]$ of $\Lambda$ at $\omega$ has been introduced in Definition~\ref{Def. Graded Kro and loc at.}. Using notations of Lemma~\ref{Lem. A/B circ P(quiver) description}, we define the morphism of DG algebras $\phi: (K\widehat{Q}/\langle I \rangle , \widehat{| \ . \ |},\widehat{d}) \rightarrow \Lambda[\omega^{-1}]$ in the following way:

        Let $\phi: Q_0 \sqcup \widehat{Q}_1 \rightarrow K\widetilde{Q}$ be given by $\phi_{|Q_0 \sqcup Q_1 \sqcup \lbrace \delta \rbrace}=id$, $\phi(\epsilon_1) = \phi(\epsilon_2) = \phi(\gamma) = 0$. It extends to an epimorphism of graded algebra $\phi: K\widehat{Q} \rightarrow K\widetilde{Q}$, and to an epimorhism of DG algebra $\phi: (K\widehat{Q}/\langle I \rangle , \widehat{| \ . \ |},\widehat{d}) \rightarrow \Lambda[\omega^{-1}]$ since $\phi(I) \subset \widetilde{I}$ and $\phi(\widehat{d}(\epsilon_1)),\phi(\widehat{d}(\epsilon_2)),\phi(\widehat{d}(\gamma)) \in \langle\widetilde{I} \rangle $.

\end{definition}

The next subsection will be devoted to the proof of the following proposition.

\begin{proposition}\label{Prop. q.e. P(phi)}
    The functor $\mathcal{P}(\phi): (\mathcal{A}/\mathcal{B})_{\circ} \rightarrow \mathcal{P}(\Lambda[\omega^{-1}])$ induced by the morphism $\phi$ of Definition~\ref{Def. loc at gr K and phi to be q.e.} is a quasi-equivalence.
\end{proposition}

\subsection{Computation of $\mathrm{H^*}(\mathcal{A}/\mathcal{B})_{\circ}$}

\begin{lemma}\label{Lem. H^* End B}
    $H^*\mathrm{End}_{\mathcal{A}}(B) \ \simeq K[x]/(x^2)$, where x is of degree 1
\end{lemma}

\textit{Proof:} For a graded pinched gentle algebra $\Lambda = (KQ / \langle I \rangle , | \ . \ |)$ and $i\neq j \in Q_0$, a basis of $e_j \Lambda e_i \simeq e_j \Lambda^g e_i$ is given by paths in $(Q^g,I^g)$. Let us give a basis of $e_1 \Lambda e_2$ in case it is non-empty. 

\vspace{3mm}

First there can be no path from 2 to 1 of the form $\rho' \alpha \rho$ or $\rho' \beta \rho$. Indeed, if $\rho' \alpha \rho \in e_1 \Lambda e_2$ then the acyclicity of $(\alpha, \beta)$ ensures that $\rho \alpha \in I^g$ and $\alpha \rho' \in I^g$, but then $\rho \beta \notin I^g$ and $\beta \rho' \notin I^g$, and thus $\rho' \alpha \rho \beta$ is an oriented cycle, a contradiction.

Now let $\rho$  be a non-zero path in $e_1 \Lambda e_2$ not passing through $\alpha$ or $\beta$. We can suppose without loss of generality that $\rho \alpha \notin I^g$. The acyclicity of $(\alpha,\beta)$ implies that $\alpha \rho \in I^g$ and thus $\beta \rho \notin I^g$. 

There can be no other path in $e_1\Lambda e_2$ different from $\rho$. Indeed suppose that $\rho'$ is such a path. Then $\rho' \beta \in I^g$, otherwise $\rho' \beta \rho \in e_1 \Lambda e_2$ is a contradiction. Since in a gentle bound quiver, every path $\gamma$ belongs to a unique maximal path $\hat{\gamma}$, we have $\hat{\rho} = \hat{\alpha} = \hat{\rho'}$. But then we can suppose that $\rho = \rho'' \rho'$, and thus $\rho = \rho''' \alpha \rho'$ or $\rho = \rho''' \beta \rho'$, a contradiction. 

In conclusion, if $e_1 \Lambda e_2 \neq \emptyset$ then $e_1 \Lambda e_2 = \langle \rho  \rangle $ where $\rho$ is a path which we can suppose verify $\rho \alpha, \beta \rho \notin I$. For the rest this proof, let $\rho = 0$ if $e_1 \Lambda e_2$ is zero.

\vspace{3mm}

We can deduce:

\vspace{-4mm}

\begin{equation} \label{Equation. e_1,2 Lambda e_1,2}
        e_1 \Lambda e_1 = \langle e_1, \rho \alpha  \rangle , \ \ e_2 \Lambda e_1 = \langle \alpha, \beta, \beta \rho \alpha  \rangle , \ \ e_1 \Lambda e_2 = \langle \rho  \rangle , \ \ e_2 \Lambda e_2 = \langle e_2, \beta\rho  \rangle ,
\end{equation}

which gives the basis $\mathrm{End}_{\mathcal{A}}(B) = \langle [e_1], [e_2], [\alpha], [\beta], [\rho], [\rho \alpha], [\beta\rho], [\beta \rho \alpha]  \rangle $,
 with 
 
 $|[e_1]|=|[e_2]| = 0,|[\alpha]|=|[\beta]|=1, |[\rho]| = |\rho| + q - 1$.

\vspace{3mm}

Dropping without ambiguity the bracket notation, the differential  is given by:

\vspace{-0.5cm}
\begin{align*}
d(e_1) = \omega = - d(e_2), \ \ & \hspace{1cm} \ \ \ \ d(\rho) = \omega \rho - (-1)^{|[\rho]|} \rho \omega = \mu \beta \rho -(-1)^{|[\rho]|} \rho \alpha, \\
d(\alpha) = d(\beta) = 0, \ \ \ \ \ & \hspace{1cm}  d(\beta \rho) = (-1)^{|[\rho]|} \beta \rho \alpha, \ \ d(\rho \alpha) =  \mu \beta \rho \alpha, \ \ d(\beta \rho \alpha) = 0,
\end{align*}

which gives the decomposition as complexes

\begin{center}
    $\mathrm{End}_{\mathcal{A}}(B) \ \simeq \ \langle [e_1], [e_2], [\alpha], [\beta]  \rangle  \bigoplus \langle [\rho], [\rho \alpha], [\beta\rho], [\beta \rho \alpha]  \rangle  \ =: \ E_0 \bigoplus E_1$.
\end{center}

Now $E_1$ is acyclic since $\mathrm{Ker} \ d_{|E_1} = \langle \rho \alpha - (-1)^{|[\rho]|} \mu \beta \rho  \rangle  \bigoplus \langle  \beta \rho \alpha  \rangle  = \mathrm{Im} \ d_{|E_1}$, and 

\vspace{-2mm}

\begin{equation}\label{Eq. H^*(B) basis}
    H^*\mathrm{End}_{\mathcal{A}}(B) \ \simeq \ H^*E_0 \ \simeq \  \langle id = e_1 + e_2  \rangle ^0 \bigoplus \langle \alpha, \beta \ | \ \omega  \rangle ^1.
\end{equation}

\hfill $\square$

A direct computation shows that:

\begin{lemma}\label{Lem. H^* (P_i,B,P_j) basis}
    Using the notations of (\ref{Equation. e_1,2 Lambda e_1,2}), 
    
    \begin{enumerate}
        \item [$\bullet$] $H^* \mathrm{Hom}_{\mathcal{A}}(P_1,B) \simeq \langle \alpha, \beta | \omega  \rangle $,
        \item [$\bullet$] $H^* \mathrm{Hom}_{\mathcal{A}}(P_2,B) \simeq \langle e_2  \rangle $,
        \item [$\bullet$] $H^* \mathrm{Hom}_{\mathcal{A}}(B,P_1) \simeq \langle e_1  \rangle $,
        \item [$\bullet$] $H^* \mathrm{Hom}_{\mathcal{A}}(B,P_2) \simeq \langle \alpha, \beta | \omega  \rangle $.
    \end{enumerate}
\end{lemma}

\begin{lemma} \label{lemma H^*(P_i,B)=0}
    For all $i \in Q_0 \backslash \lbrace 1, 2 \rbrace$,    $H^*\mathrm{Hom}_{\mathcal{A}} (P_i, B) \simeq 0 \simeq H^*\mathrm{Hom}_{\mathcal{A}}(B,P_i)$.
\end{lemma}

\textit{Proof:} Let $e_1 \Lambda e_i \simeq e_1 \Lambda^g e_i = \langle \rho_1, ..., \rho_r  \rangle $ and $e_2 \Lambda e_i \simeq e_2 \Lambda^g e_i = \langle \gamma_1, ..., \gamma_s  \rangle $ be two bases where each $\rho_k, \gamma_l$ is a path in $(Q^g,I^g)$. We use the same notation to designate the induced basis of $\mathrm{Hom}_{\mathcal{A}}(P_i,B)$.

Recall that $B$ is of the form $(P_2 [|\alpha|] \bigoplus P_1 [1], \partial =
        \begin{pmatrix}
        0 & \alpha + \mu \beta \\
        0 & 0
        \end{pmatrix} )$ for some $\mu \in K^*$, and thus for $f$ in $\mathrm{Hom}_{\mathcal{A}}(P_i,B)$ the differential is given by $df = (\alpha + \mu \beta)f$ (see Notations~\ref{Not. DG recollections}).

\vspace{2mm}

Since $i \neq 1$, one can choose an ordering of the $\rho_l$'s and $r' \in \lbrace 0, ..., r \rbrace$ such that $\forall l \in \lbrace 1, ..., r' \rbrace, \ \beta \rho_l \in I$ and $\forall l \in \lbrace r' + 1, ..., r \rbrace, \ \alpha \rho_l \in I$. So for $(\lambda_l),(\mu_h) \in K^{r+s}$,

\vspace{-4mm}

\begin{align*}
    d(\sum\limits_{l=1}^r \lambda_l \rho_l + \sum\limits_{h=1}^s \mu_h \gamma_h) = 0 \ & \Leftrightarrow \ (\alpha + \mu \beta)(\sum\limits_{l=1}^{r'} \lambda_l \rho_l + \sum\limits_{r' + 1}^r \lambda_l \rho_l) + (\alpha + \mu \beta)(\sum\limits_{h=1}^s \mu_h \gamma_h) = 0 \\
    & \Leftrightarrow \ \sum\limits_{l=1}^{r'} \lambda_l \alpha \rho_l + \sum\limits_{r' + 1}^r \lambda_l \mu \beta \rho_l = 0 \\
    & \Leftrightarrow \ \forall l \in \lbrace 1, ..., r \rbrace, \ \lambda_l = 0,
\end{align*}

where the last equivalence comes from the fact that in a gentle algebra distinct paths correspond to distinct elements, and thus the $\alpha\rho_l$'s and the $\beta\rho_l$'s are pairwise distinct. We  deduce that $\mathrm{Ker} \ d = \langle \gamma_1, ..., \gamma_s  \rangle $.

\vspace{2mm}

Now since $i \neq 2$, one can choose an ordering of the $\gamma_h$'s and $s' \in \lbrace 0, ..., s \rbrace$ such that $\forall h \in \lbrace 1, ..., s' \rbrace, \exists l_h \in \lbrace 1, ..., r \rbrace, \ \gamma_h = \alpha \rho_{l_h}$ and $\forall h \in \lbrace s' + 1, ..., s \rbrace,  \exists l_h \in \lbrace 1, ..., r \rbrace, \ \gamma_h = \beta \rho_{l_h}$. Thus for all $(\mu_h) \in K^s$,

$\sum\limits_{h=1}^s \mu_h \gamma_h = \sum\limits_{h=1}^{s'} \mu_h \alpha \rho_{l_h} + \sum\limits_{h=s' + 1}^s \frac{\mu_h}{\mu} \mu \beta \rho_{l_h} = d(\sum\limits_{h=1}^{s'} \mu_h \rho_{l_h} + \sum\limits_{h = s' + 1}^s \frac{\mu_h}{\mu} \rho_{l_h})$.

\vspace{2mm}

We can deduce $\mathrm{Ker} \ d \simeq \mathrm{Im} \ d$ and $H^*\mathrm{Hom}_{\mathcal{A}}(P_i,B) = 0$. The other case is dual.

\hfill $\square$

\begin{corollary}\label{Cor. i,j neq 1,2 Hom A/B}
    If i or j does not belong to $\lbrace 1,2 \rbrace$, then $H^*\mathrm{Hom}_{\mathcal{A}/\mathcal{B}}(P_i,P_j) \simeq \mathrm{Hom}_{\Lambda}(P_i,P_j)$.
\end{corollary}

\textit{Proof:} Let $(i,j)$ be such a pair. By Equation~\ref{Eq. H* DG quot via Ss}, $H^{k} \mathrm{Hom}_{\mathcal{A}/\mathcal{B}}(P_i,P_j) \simeq \bigoplus\limits_{l\in\mathbb{Z}} E_{\infty}^{k-l,l}$.

Let $q \in \mathbb{Z}$. Since $H^*\mathrm{Hom}_{\mathcal{A}} (P_i, B) \simeq 0$ or $H^*\mathrm{Hom}_{\mathcal{A}}(B,P_j) \simeq 0$ by Lemma~\ref{lemma H^*(P_i,B)=0}, Proposition~\ref{prop E^1} shows that $\forall p > 0$, $E_1^{p,q} \simeq 0$. Thus $\forall p \neq 0$, $E_{\infty}^{p,q} \simeq 0$. 

Moreover since for all $r > 1$ the codomain of $d_r^{0,q}$ and the domain of $d_r^{r, q-r-1}$ are zero, we have $E_{\infty}^{0,q} \simeq E_1^{0,q} \simeq \mathrm{Hom}_{\Lambda}(P_i,P_j)^q$. 

\hfill $\square$

\begin{proposition}\label{Prop. E_infty computation}
    Let $(i,j) \in \lbrace (1,1),(1,2),(2,1),(2,2) \rbrace$ and $E^{*,*}$ the spectral sequence on $\mathrm{Hom}_{\mathcal{A}/\mathcal{B}}(P_i,P_j)$.

     Let $H^* \mathrm{Hom}_{\mathcal{A}}(P_i,B) =: \langle t  \rangle $, $H^*\mathrm{Hom}_{\mathcal{A}}(B,P_j) =: \langle z  \rangle $ and $H^*\mathrm{End}_{\mathcal{A}}(B) =: \langle id_B  \rangle  \bigoplus \langle y  \rangle $ be the basis given in Lemma~\ref{Lem. H^* (P_i,B,P_j) basis} and in (\ref{Eq. H^*(B) basis}), and let $a:= |z| + |t| -1$
     
    \begin{enumerate}
        \item [$\bullet$] $\forall q \in \mathbb{Z}$, $E_{\infty}^{0,q} \simeq \mathrm{Hom}_{\Lambda}(P_i,P_j)^q$,
        \item [$\bullet$] $\forall p > 0$, $E_{\infty}^{p, a-p} \simeq \langle z \epsilon (y \epsilon)^{p-1} t  \rangle $,
        \item[$\bullet$] $\forall p \neq 0, \forall q \neq a - p$, $E_{\infty}^{p,q} \simeq 0$.
    \end{enumerate}
\end{proposition}

\textit{Proof:} We have the following relations: $zy = yt = y^2 = zt = 0$. Let $p>0$. For $q \in \mathbb{Z}$, since by Proposition~\ref{prop E^1}

\vspace{-3mm}

\begin{equation*}
E_1^{p,q} \simeq \bigoplus\limits_{\substack{k_1 + ... + k_{p+1} \\ = \ 2p + q}} H^{k_{p+1}}(\mathrm{Hom}_{\mathcal{A}}(B,P_j)) \otimes K[1] \otimes H^{k_p} (\mathrm{End}_{\mathcal{A}}(B)) \otimes ... \otimes K[1] \otimes H^{k_1} (\mathrm{Hom}_{\mathcal{A}}(P_i,B)),
\end{equation*}

we have 

\vspace{-4mm}
\begin{align*}
E_1^{p,q} \neq 0 & \Leftrightarrow \exists (k_i) \in \lbrace 0, 1 \rbrace^{p-1}, \ |z| + \sum\limits_{i} k_i + |t| = 2p+q \\
& \Leftrightarrow 0 \leq 2p + q - |z| - |t| \leq p-1 \\
& \Leftrightarrow a + 1 -2p \leq q \leq a - p.
\end{align*}

For such $p$ and $q$, a basis of $E_1^{p,q}$ is given by 

$E_1^{p,q} = \langle \lbrace z \epsilon^{n_0} \Pi_{l=1}^r (y \epsilon^{n_l})t \ | \ r \geq 0, \ \forall l \in \lbrace 0 , ..., r \rbrace. \ n_l \geq 1, \ \sum\limits_{l=0}^r n_l = p \ \mathrm{and} \ |z| + r -p + |t| = p + q \rbrace  \rangle $,

where $\epsilon$ is a generator of $K[1]$. Moreover,

\vspace{-3mm}
\begin{align*}
d_1^{p,q}(z \epsilon^{n_0} \Pi_{l=1}^r (y \epsilon^{n_l})t) & = \sum\limits_{i = 0}^r (-1)^{| z (\Pi_{l=1}^i \epsilon^{n_{l-1}} y) |} z(\Pi_{l=1}^i \epsilon^{n_{l-1}} y) d(\epsilon^{n_i}) (\Pi_{l=i+1}^r y \epsilon^{n_l})t \\
& = \sum\limits_{\substack{i=0 \\ n_i \neq 1 \ \mathrm{odd}}}^r (-1)^{|z (\Pi_{l=1}^i \epsilon^{n_{l-1}} y)|} z (\Pi_{l=1}^i \epsilon^{n_{l-1}} y) \epsilon^{n_i - 1} (\Pi_{l=i+1}^r y \epsilon^{n_l})t,
\end{align*}

since $\forall n \geq 1. \ d(\epsilon^n) = \sum\limits_{i=1}^n (-1)^{|\epsilon^{i-1}|} \epsilon^{n-1} = 0$ if $n$ even and $\epsilon^{n-1}$ otherwise.

\vspace{5mm}

We first compute the second page of the spectral sequence.

\vspace{3mm}

For $q \in \mathbb{Z}$, $d_1^{0,q} = 0$ and $E_1^{1,q} \neq 0 \Leftrightarrow q = a-1$. Since $d_1^{1,a-1}(z \epsilon t) = (-1)^{|z|}zt = 0$, $E_2^{0,q} \simeq E_1^{0,q}.$

\vspace{3mm}

For all $p > 0$, $d_1^{p,a-p}(z \epsilon (y \epsilon)^{p-1} t) = 0$, and for all $\rho = z \epsilon^{n_0} \Pi_{l=1}^r (y \epsilon^{n_l})t \in E_1^{p+1,a-p-2}$, $d_1^{p+1,a-p-2}(\rho) = 0$ since $n_l = 1$ or $2$. Thus $E_2^{p,a-p} \simeq E_1^{p,a-p}$.

\vspace{3mm}

We regroup the remaining cases by complexes. For each $b>0$, let $C_b^{\bullet}$ be defined by $\forall i\geq a, \ C_b^i = 0$, $\forall i<a, \ C_b^i = E_1^{b+(a-i),a-b-2(a-i)}$ and $\forall i<a, \ d_C^i = d_1^{b+(a-i), a-b-2(a-i)}$. In this way, $\forall p>0$ and $\forall q \in \lbrace a +1 -2p, ..., a-p-1 \rbrace$, $E_1^{p,q}$ is equal to $C_{2p + q -a}^{p+q}$. Let us show that each $C_b^{\bullet}$ is acyclic by defining morphisms $\Psi: C_b^{\bullet} \rightarrow C_b^{\bullet} [-1]$ satisfying $d_{C_b^{\bullet}} \circ \Psi + \Psi \circ d_{C_b^{\bullet}} = id_{C_b^{\bullet}}$.

For $i < a$ and $\rho = z \epsilon^{n_0} \Pi_{l=1}^r (y \epsilon^{n_l})t$ a basis element of $C_b^i$, let $X = \lbrace l \in \lbrace 0,..., r \rbrace | \ n_l \ \mathrm{is \ odd \ and \ different \ from} \ 1\rbrace$, $Y = \lbrace l \in \lbrace 0, ..., r \rbrace | \ n_l \ \mathrm{is \ even} \rbrace$ and $\lambda = \frac{1}{|X| + |Y|}$.

Define $\Psi$ on paths by $\Psi(\rho) = \lambda \sum\limits_{i \in Y} (-1)^{|z (\Pi_{l=1}^i \epsilon^{n_{l-1}}y)|}z (\Pi_{l=1}^i \epsilon^{n_{l-1}}y) \epsilon^{n_i +1} (\Pi_{l=i+1}^r y \epsilon^{n_l})t$.

Then $\rho$ appears with coefficient $\lambda|Y|$ in $d \circ \Psi (\rho)$ and with coefficient $\lambda |X|$ in $\Psi(\rho) \circ d$. For any other path $\rho' \neq \rho$ appearing with coefficient $\mu$ in $d \circ \Psi (\rho)$, it appears with coefficient $-\mu$ in $\Psi \circ d (\rho)$ and vice-versa. Thus  $d \circ \Psi + \Psi \circ d = \rho$ holds.

Finally the values of $E_{\infty}$ correspond to the second page since $\forall p,q \in \mathbb{Z}, \forall r \geq 2, d_r^{p,q} = 0$. 

\hfill $\square$

\textit{Proof of Proposition~\ref{Prop. q.e. P(phi)}:} Since $\mathcal{P}(\Lambda[\omega^{-1}])$ has zero differential, it suffices to show that $H^* \mathcal{P}(\phi): H^* \mathrm{Hom}_{(\mathcal{A}/\mathcal{B})_{\circ}}(P_i,P_j) \rightarrow \mathrm{Hom}_{\mathcal{P}(\Lambda[\omega^{-1}])}(P_i,P_j)$ is an isomorphism for all $i,j$.

\vspace{1mm}

If $i$ or $j$ is different from 1 or 2, the relations in $\Lambda[\omega^{-1}]$ are such that $e_j \Lambda[\omega^{-1}] e_i \simeq e_j \Lambda e_i$. Thus $H^* \mathcal{P}(\phi)$ is an isomorphism by Corollary~\ref{Cor. i,j neq 1,2 Hom A/B}.

\vspace{1mm}

For $(i,j) \in \lbrace (1,1),(1,2),(2,1),(2,2) \rbrace$, Equation~\ref{Eq. H* DG quot via Ss} and Proposition~\ref{Prop. E_infty computation}
 give the isomorphism:

\begin{center}
    $H^*\mathrm{Hom}_{\mathcal{A}/\mathcal{B}}(P_i,P_j) \simeq \mathrm{Hom}_{\Lambda}(P_i,P_j) \bigoplus \langle z \epsilon (y \epsilon)^n t \ | \ n \geq 0  \rangle . $
\end{center}

\vspace{1mm}

Let us consider the case $(i,j) = (1,1)$. On one hand, the relations in $\Lambda[\omega^{-1}]$ allow to choose the following basis: $\mathrm{Hom}_{\mathcal{P}(\Lambda[\omega^{-1}])} (P_1,P_1) = \mathrm{Hom}_{\Lambda}(P_1,P_1) \bigoplus \langle \delta (\alpha \delta)^n \alpha \ | \ n \geq 0  \rangle $. On the other hand, if we choose the following representatives: $t = \alpha, y = \alpha, z = e_1$, the equivalence of Lemma~\ref{Lem. A/B circ P(quiver) description} gives the rewriting $z \epsilon (y \epsilon)^n t = e_1 \epsilon (e_2 \alpha e_1 \epsilon)^n e_2 \alpha = \delta (\alpha \delta)^n \alpha$, and $H^* \mathcal{P}(\phi)$ is an isomorphism.

\vspace{1mm}

Similarly for the other cases we found:

\begin{enumerate}
    \item[] $H^*\mathrm{Hom}_{\mathcal{A}/\mathcal{B}}(P_2,P_2) = \mathrm{Hom}_{\Lambda}(P_2,P_2) \bigoplus \langle \alpha \delta (\alpha \delta)^n | \ n \geq 0  \rangle $,
    \item[] $H^*\mathrm{Hom}_{\mathcal{A}/\mathcal{B}}(P_1,P_2) = \mathrm{Hom}_{\Lambda}(P_1,P_2) \bigoplus \langle \alpha \delta (\alpha \delta)^n \alpha | \ n \geq 0  \rangle $,
    \item[] $H^*\mathrm{Hom}_{\mathcal{A}/\mathcal{B}}(P_2,P_1) = \mathrm{Hom}_{\Lambda}(P_2,P_1) \bigoplus \langle \delta (\alpha \delta)^n | \ n \geq 0  \rangle $,
\end{enumerate}

which are isomorphic under $H^* \mathcal{P}(\phi)$ to the corresponding spaces in $\mathcal{P}(\Lambda[\omega^{-1}])$.

\hfill $\square$

\section{\textbf{Proof of the main results}}\label{Sec. pf main res}

We now prove Theorem~\ref{Prop. Pinched Main quotient} and~\ref{Th. Pinched Main recollement}. For this we will need two lemmas.

\begin{lemma}\label{Lem. subalg quiver}
    Let $\Lambda$ be the quotient of a path algebra $KQ$ by a finitely generated ideal $I = \langle \rho_1,..., \rho_r  \rangle $. Let $\overline{Q}_0 \subseteq Q_0$ be a subset of vertices and $e = \sum\limits_{i \in \overline{Q}_0} e_i$ the associated idempotent. 

    There is an isomorphism $[[-]]: e KQ e \rightarrow K\overline{Q}$, where $\overline{Q}$ is the quiver whose set of vertices is $\overline{Q}_0$, and whose arrows from $i$ to $j$ are symbols $[\alpha_r ... \alpha_1]$, where $\alpha_r ... \alpha_1$ is a path from $i$ to  $j$ in $Q$ which passes through no other vertices of $\overline{Q}_0$.
    \vspace{1mm}
    
    It induces an isomorphism $e \Lambda e \simeq K \overline{Q} / \overline{I}$, where $\overline{I} = \langle [[\delta \rho_i \gamma]] \ | \ \forall i \in \lbrace 1,...,r \rbrace \ \mathrm{and} \ \delta \in T, \ \gamma \in S  \rangle $, with:

    \begin{enumerate}
        \item[-] $T = \lbrace \alpha_r ... \alpha_1 \ | \ r \geq 1, \ \forall i. \alpha_i \in Q_1 \ \mathrm{and} \ \tau(\alpha_r) \in \overline{Q}_0 \ \mathrm{and} \ \forall i. \sigma(\alpha_i) \in Q_0 \backslash \overline{Q}_0 \rbrace \bigsqcup \lbrace e_i \ | \ i \in \overline{Q}_0 \rbrace$, and dually
        \item[-] $S = \lbrace \alpha_r ... \alpha_1 \ | \ r \geq 1, \ \forall i. \alpha_i \in Q_1 \ \mathrm{and} \ \sigma(\alpha_1) \in \overline{Q}_0 \ \mathrm{and} \ \forall i. \tau(\alpha_i) \in Q_0 \backslash \overline{Q}_0 \rbrace \bigsqcup \lbrace e_i \ | \ i \in \overline{Q}_0 \rbrace$.
    \end{enumerate}
    
\end{lemma}

Note that $\overline{Q}$ can have an infinite set of arrows, and that $\overline{I}$ is not necessarily admissible.

\begin{lemma}\label{Lem. eq D(loc at K) D(pinch at K)}
    Let $\Lambda = (KQ/\langle I  \rangle , | \ . \ |)$ be a graded pinched gentle algebra, $(\alpha, \beta)$ be a graded Kronecker in $Q$ as in Definition~\ref{Def. Pinched Graded Kro and loc at.}, and $\omega = \alpha + \mu \beta$ for some $\mu \in K^*$. Suppose that the characteristic of $K$ is different from 2. There are equivalences:

    \begin{center}
        $\mathcal{D}(\Lambda[\omega^{-1}]) \simeq \mathcal{D}(\Lambda_{(\alpha, \beta)})$ \hspace{2mm} and \hspace{2mm} $per(\Lambda[\omega^{-1}]) \simeq per(\Lambda_{(\alpha, \beta)})$.
    \end{center}
\end{lemma}

\textit{Proof:} Since $P_1$ and $P_2$ become isomorphic in $\Lambda[\omega^{-1}] =: K\widetilde{Q} / \langle \widetilde{I}  \rangle $, $\mathcal{D}(\Lambda[\omega^{-1}]) \simeq \mathcal{D}(e \Lambda[\omega^{-1}] e)$ and $per(\Lambda[\omega^{-1}]) \simeq per(e \Lambda[\omega^{-1}] e)$ where $e = \sum\limits_{l \in Q_0 \backslash \lbrace 2 \rbrace} e_l$.

One can compute $e \Lambda[\omega^{-1}] e =: K\overline{Q}/ \overline{I} $ using Lemma~\ref{Lem. subalg quiver}. First we have:

\vspace{2mm}

$\overline{Q}_1 = \lbrace [\rho] \ | \ \rho \in \widetilde{Q}_1 \ \mathrm{and} \ \widetilde{\sigma}(\rho) \neq 2 \neq \widetilde{\tau}(\rho) \rbrace \bigsqcup \lbrace [\beta^+ \beta],[\alpha^+ \alpha],[\delta\alpha],[\delta\beta] \rbrace$,

\vspace{2mm}

and Lemma~\ref{Lem. subalg quiver} gives the following description of $\overline{I}$: 

\vspace{-4mm}

\begin{align*}
    \overline{I} = & \langle [[\delta \rho_i \gamma]] \ | \rho_i \in \widetilde{I} = I^g \sqcup I^p \sqcup \lbrace \delta \omega - e_1, \omega \delta - e_2 \rbrace \ \mathrm{and} \ \delta \in T, \ \gamma \in S  \rangle  \\
    = & \langle \lbrace [\beta_s ... \beta_1] [\alpha_r ... \alpha_1] \ | \ [\beta_s ... \beta_1], [\alpha_r ... \alpha_1] \in \overline{Q}_1 \ \mathrm{and} \ \beta_1 \alpha_r \in I^g \rbrace \\
    \sqcup & \lbrace [\beta_s ... \beta_1]([\gamma_v] +e_v) \ | \ [\beta_s ... \beta_1] \in \overline{Q}_1 \ \mathrm{and} \ \beta_1(\gamma_v + e_v) \in I^p \rbrace \\
    \sqcup & \lbrace ([\gamma_v] + e_v)[\beta_s ... \beta_1] \ | \ [\beta_s ... \beta_1] \in \overline{Q}_1 \ \mathrm{and} \ (\gamma_v + e_v)\beta_s \in I^p \rbrace \\
    \sqcup & \lbrace [\alpha_r ... \alpha_1]([\gamma_v] - e_v) \ | \ [\alpha_r ... \alpha_1] \in \overline{Q}_1 \ \mathrm{and} \ \alpha_1(\gamma_v - e_v) \in I^p \rbrace \\
    \sqcup & \lbrace ([\gamma_v] - e_v)[\alpha_r ... \alpha_1] \ | \ [\alpha_r ... \alpha_1] \in \overline{Q}_1 \ \mathrm{and} \ (\gamma_v - e_v)\alpha_s \in I^p \rbrace \\
    \sqcup & \lbrace [[ \delta \omega - e_1 ]] = [\delta \alpha] + \mu [\delta \beta] - e_1, \\
    & \ [[ \alpha^+ (\omega \delta - e_2) \alpha]] = [\alpha^+ \alpha][\delta \alpha] - [\alpha^+ \alpha], \ [[ \alpha^+ (\omega \delta - e_2) \beta]] = [\alpha^+ \alpha][\delta \beta],\\
    & \ [[ \beta^+ (\omega \delta - e_2) \alpha]] = \mu[\beta^+ \beta][\delta \alpha], \ [[ \beta^+ (\omega \delta - e_2) \beta]] = \mu[\beta^+ \beta][\delta\beta] - [\beta^+\beta] \rbrace \rangle.    
\end{align*}

Let $\gamma:= [\delta\alpha] - \mu[\delta\beta]$. We have $2[\delta\alpha] = e_1 + [\delta \alpha] - \mu [\delta \beta] = e_1 + \gamma$ and $2\mu[\delta\beta] = e_1 -([\delta\alpha] - \mu [\delta\beta]) = e_1 - \gamma$. Renaming $[\beta^+ \beta] =: \beta^+$, $[\alpha^+ \alpha] =: \alpha^+$, dropping the bracket notation and rewriting the relations with $\gamma's$ gives the desired description.

\hfill $\square$

\textit{Proof of Theorem~\ref{Prop. Pinched Main quotient}:} By Theorems \ref{Th Dri Quo}, \ref{Th. BK thick} and Example~\ref{Ex. P Lambda tr simeq per}, $per(\Lambda) / \langle B  \rangle  \simeq \mathcal{A}^{tr} / \mathcal{B}^{tr} \simeq (\mathcal{A} / \mathcal{B} ) ^{tr}$, and by Proposition~\ref{Prop. q.e. P(phi)}, $(\mathcal{A} / \mathcal{B} ) ^{tr} \simeq (\mathcal{A}/\mathcal{B})_{\circ}^{tr} \simeq \mathcal{P}(\Lambda[\omega^{-1}])^{tr} \simeq per(\Lambda[\omega^{-1}])$. We conclude using $per(\Lambda)/\langle B \rangle \simeq per(\Lambda)/thick(B)$.

The second part of the proposition is given by Lemma~\ref{Lem. eq D(loc at K) D(pinch at K)}.

\hfill $\square$

\textit{Proof of Theorem~\ref{Th. Pinched Main recollement}:} Let $\mathcal{A}$ and $\mathcal{B}$ be as in the Setting~\ref{Set. mc(A)}. Theorem 3.1 of \cite{G22} gives the recollement:

\begin{equation*}
    \begin{tikzcd}
    \mathcal{D}(\mathcal{A}/\mathcal{B}) \arrow[r] & \mathcal{D}(\mathcal{A}) \arrow[l, shift left=0.75ex] \arrow[l, shift right=0.75ex] \arrow[r] & \mathcal{D}(\mathcal{B}) \arrow[l, shift left=0.75ex] \arrow[l, shift right=0.75ex]
    \end{tikzcd}
\end{equation*}

First $\mathcal{D}(\mathcal{A}) \simeq \mathcal{D}(\mathcal{P}(\Lambda)) \simeq \mathcal{D}(\Lambda)$. Then one can show that $\mathrm{End}_{\mathcal{P}(\Lambda)^{pre-tr}}(B)$ is a formal DG algebra, and thus $\mathcal{D}(\mathcal{B}) \simeq \mathcal{D}(\mathrm{End}_{\mathcal{P}(\Lambda)^{pre-tr}}(B)) \simeq \mathcal{D}(H^* \mathrm{End}_{\mathcal{P}(\Lambda)^{pre-tr}}(B)) \simeq \mathcal{D}(K[x]/(x^2))$, by Lemma~\ref{Lem. H^* End B}.

Moreover $\mathcal{D}(\mathcal{A}/\mathcal{B}) \simeq \mathcal{D}((\mathcal{A}/\mathcal{B})_{\circ}) \simeq \mathcal{D}(\Lambda[\omega^{-1}])$, where the last equivalence is given by Proposition~\ref{Prop. q.e. P(phi)}. 

The second part of the theorem is given by Lemma~\ref{Lem. eq D(loc at K) D(pinch at K)}.

\hfill $\square$


\end{document}